\documentclass[a4paper,12pt]{article}
\usepackage{amsthm,amsmath}
\usepackage{amssymb}
\usepackage[all]{xypic}
\usepackage{enumerate}
\usepackage{a4wide}
\usepackage{hyperref}
\usepackage{latexsym}
\usepackage{enumerate}

\newtheorem{theorem}{Theorem}[section]
\newtheorem{proposition}[theorem]{Proposition}
\newtheorem{corollary}[theorem]{Corollary}
\newtheorem{lemma}[theorem]{Lemma}
\newtheorem*{theorem*}{Theorem}
\newtheorem*{proposition*}{Proposition}
\newtheorem*{corollary*}{Corollary}
\newtheorem*{lemma*}{Lemma}
\theoremstyle{definition}
\newtheorem{definition}[theorem]{Definition}

\newtheorem{example}[theorem]{Example}
\newtheorem{remark}[theorem]{Remark}
\newtheorem*{remark*}{Remark}

\newtheorem*{definition*}{Definition}

\newcommand{\cat}[1]{\mathcal{#1}}
\newcommand{\coring}[1]{\mathfrak{#1}}
\newcommand{\tensor}[1]{\otimes_{#1}}

\newcommand{\rcomod}[1]{\mathcal{M}^{#1}}
\newcommand{\rmod}[1]{\mathcal{M}_{#1}}
\newcommand{\rmd}[1]{\widetilde{\mathcal{M}}_{#1}}
\newcommand{\lmod}[1]{{}_{#1}\mathcal{M}}

\newcommand{\cotensor}[1]{\square_{#1}}
\newcommand{\lcomod}[1]{{}^{#1}\mathcal{M}}

\renewcommand{\hom}[3]{\mathrm{Hom}_{#1}(#2,#3)}

\newcommand{\rend}[2]{\mathrm{End}({#2}_{#1})}
\newcommand{\lend}[2]{\mathrm{End}({}_{#1}#2)}

\newcommand{\rcomatrix}[2]{#2^* \tensor{#1} #2}
\newcommand{\bicomod}[2]{{}^{#1} \mathcal{M}^{#2}}

\newcommand{\dostensor}[3]{#1 \tensor{#2} #3}
\newcommand{\trestensor}[5]{#1 \tensor{#2} #3 \tensor{#4} #5}
\newcommand{\fourtensor}[7]{#1 \tensor{#2} #3 \tensor{#4} #5 \tensor{#6} #7}
\newcommand{\abrir}[1]{e_{#1}\tensor{A}\varphi_{#1}}

\def\*C{{^*\mathcal C}}
\def\*Dd{{^*\mathcal D}}

\newcommand{\Hom}{{\rm Hom}}

\newcommand{\Ker}{{\rm Ker}\,}

\newcommand{\can}{{\sf can}}

\def\ZZ{{\mathbb Z}}
\def\QQ{{\mathbb Q}}

\newcommand{\Cc}{\mathcal{C}}

\newcommand{\Mm}{\mathcal{M}}

\newcommand{\Tt}{\mathcal{T}}

\def\*C{{}^*\hspace*{-1pt}{\Cc}}

\begin{document}
\title{Comatrix Corings and Galois Comodules over firm rings}
\author{J. G\'omez-Torrecillas \\ \normalsize Departamento de \'{A}lgebra \\ \normalsize Universidad
de Granada\\ \normalsize E-18071 Granada, Spain \\  \normalsize
e-mail: \textsf{gomezj@ugr.es} \and J. Vercruysse \\ \normalsize
Faculty of Engineering \\ \normalsize Vrije Universiteit Brussel, VUB \\ \normalsize B-1015 Brussels, Belgium \\
\normalsize e-mail: \textsf{jvercruy@vub.ac.be}}

\date{}

\maketitle

\section*{Introduction}
Galois corings with a group-like element \cite{Brzezinski:2002}
provide a neat framework to understand the analogies between
several theories like the Faithfully Flat Descent for
(noncommutative) ring extensions \cite{Nuss:1997}, Hopf-Galois
algebra extensions \cite{Schneider:1990}, or noncommutative Galois
algebra extensions \cite{Kanzaki:1964,DeMeyer:1965}. A Galois
coring is isomorphic in a canonical way to the Sweedler's
canonical coring $A \tensor{B} A$ associated to a ring extension
$B \subseteq A$, where $B$ is a subring of ``coinvariants''.
Canonical means here that the group-like of the Galois coring
corresponds to the group-like $1 \tensor{B} 1$ of $A \tensor{B}
A$, where $1$ denotes the unit element of the ring $A$. Canonical
corings were used in \cite{Sweedler:1975} to formulate a predual
to the Jacobson-Bourbaki theorem for division rings.

Comatrix corings were introduced in \cite{ElKaoutit/Gomez:2003} to
make out the structure of cosemisimple corings over an arbitrary
ground ring $A$ \cite[Theorem 4.4]{ElKaoutit/Gomez:2003}. The
construction has been used in different ``Galois-type'' contexts,
where no group-like element is available, as the characterization
of corings having a projective generator \cite[Theorem
3.2]{ElKaoutit/Gomez:2003}; the tightly related non-commutative
descent for modules \cite[Theorem 3.2,Theorem
3.10]{ElKaoutit/Gomez:2003}, \cite[Theorem 3.7, Theorem 3.8,
Theorem 3.9]{Caenepeel/DeGroot/Vercruysse:unp2005},
\cite[18.27]{Brzezinski/Wisbauer:2003}; the formulation of a
predual to the Jacobson-Bourbaki correspondence for simple
artinian rings \cite{Cuadra/Gomez:unp2005}; or the construction of
a Brauer Group using corings \cite{Caenepeel/Femic:unp2005}. The
original definition of a comatrix coring was built on a bimodule
$\Sigma$ over unital rings $B$ and $A$ such that $\Sigma_A$ is
finitely generated and projective. This finiteness condition seems
to be, at a first glance, essential to define comatrix corings and
to introduce the notion of a Galois coring (see
\cite{Brzezinski/Gomez:2003}). Nevertheless, it was soon
discovered \cite{ElKaoutit/Gomez:2004} that the concepts and
results from \cite[\S 2, \S 3]{ElKaoutit/Gomez:2003}, including
comatrix and Galois corings (or Galois comodules, as in
\cite{Brzezinski/Wisbauer:2003,Brzezinski:unp2004}) and the
Descent Theorem, can be developed for a functor from a small
category to a category of finitely generated and projective
modules. These ``infinite'' comatrix corings have been recently
constructed in a more general functorial setting in
\cite{DeGroot:2005}.

In this paper we recover the idea of focusing on a single
bimodule, relaxing the condition on the existence of a unity in
the rings. This allows the definition of a concept of comatrix
coring which embodies all former constructions and, what is more
interesting, leads to the formulation of  a notion of Galois
coring and the statement of a Faithfully Flat Descent Theorem that
generalize the aforementioned previous versions.

The paper is organized as follows. After fixing some basic
notations and concepts in Section \ref{basics}, we state in
Section \ref{pairsadj} some adjunctions between categories of
modules which will be used later. These adjunctions, and the
construction of the comatrix coring in Section
\ref{comatrixsection}, start from a homomorphism of non unital
rings $\iota : R \rightarrow S$, where $R$ is a firm ring, and the
ring $S = \Sigma \tensor{A} \Sigma'$ is built from two bimodules
${}_B\Sigma_A$ and ${}_A\Sigma'_B$ over unital rings $A$, $B$, and
a homomophism of $A$--bimodules $\mu : \Sigma' \tensor{B} \Sigma
\rightarrow A$. The ``finite'' comatrix corings from
\cite{ElKaoutit/Gomez:2003} fit in this framework taking $B = R$
and $\Sigma' = \Sigma^* = \hom{A}{\Sigma}{A}$ since, in that case,
$S \cong \rend{A}{\Sigma}$ because $\Sigma_A$ is assumed to be
finitely generated and projective. In the present approach, this
finiteness condition is replaced by the requirement on $\Sigma$ of
being firm as a left $R$--module. The new framework requires a
different approach (including proofs) than that of
\cite{ElKaoutit/Gomez:2003} or even \cite{ElKaoutit/Gomez:2004}.
As a counterweight for this technical effort, the range of
situations covered is sensibly wider, as will be shown in sections
\ref{comparison} and \ref{rationality}. Thus, theorems
\ref{predescensoplano} and \ref{descensoplano} improve the main
results of \cite[Section 5]{ElKaoutit/Gomez:2004}, and the
equivalence between the categories of right
$\coring{C}$--comodules and of firm right modules over the ring
with local units $({}^*\coring{C})^{\rm rat}$ given in
\cite{A,CVW2} is deduced from our general results.

Section \ref{Galoiscomodules} contains the main results of the
paper. Starting from an $R-\coring{C}$--bicomodule $\Sigma$, where
$R$ is a firm ring, and $\coring{C}$ is an  $A$--coring, and
assuming that ${}_R\Sigma$ is firm, we prove that there is a
canonical map $\mathsf{can} : \rcomatrix{R}{\Sigma} \rightarrow
\coring{C}$ defined on the comatrix coring
$\rcomatrix{R}{\Sigma}$, which is a homomorphism of $A$--corings
(Proposition \ref{can1}). This canonical map allows to consider
not necessarily finitely generated Galois comodules (Definition
\ref{Galoisdef}), generalizing the former notions from
\cite{Brzezinski/Wisbauer:2003, ElKaoutit/Gomez:2003,
ElKaoutit/Gomez:2004}. This definition of a Galois comodule, far
from being merely a formal generalization, plays a central role in
the formulation of a ``Flat Descent Theory"
 for comodules,
as stated in theorems \ref{flatdescent} and
\ref{descensofielmenteplano}.

\section{Basic notation}\label{basics}

Let $K$ be a commutative ring with unity. Our additive categories
and our functors will be $K$--linear. The identity arrow at an
object $X$ of a category will be denoted by $X$.

All bimodules over $K$ are assumed to be centralized by $K$. The
term ring refers to an associative ring, although we do not
require to our rings to have a unity. All rings will be assumed to
be $K$--rings, in the sense that they are (unital) $K$--bimodules,
with a $K$--bilinear product. Rings with unity are then
$K$--algebras. Modules over rings with unity are assumed to be
unital. If $A$ is a ring with unity, then by an $A$--ring we will
understand an $A$--bimodule $R$ endowed with an associative
multiplication defined by an $A$--bimodule homomorphism $\nabla :
\dostensor{R}{A}{R} \rightarrow R$. Associativity means here that
$\nabla \circ (\dostensor{\nabla}{A}{R}) = \nabla \circ
(\dostensor{R}{A}{\nabla})$.

Dually, an $A$--coring is an $A$--bimodule $\coring{C}$ endowed
with a coassociative comultiplication defined by a homomorphism
 of $A$--bimodules $\Delta : \coring{C} \rightarrow \dostensor{\coring{C}}{A}{\coring{C}}$.
Coassociativity means here that
$(\dostensor{\Delta}{A}{\coring{C}}) \circ \Delta =
(\dostensor{\coring{C}}{A}{\Delta}) \circ \Delta$. Our corings
will always have a \emph{counity}, that is, an $A$--bimodule
homomorphism $\epsilon : \coring{C} \rightarrow A$ such that
$(\dostensor{\coring{C}}{A}{\epsilon})\circ \Delta =
(\dostensor{\epsilon}{A}{\coring{C}})\circ \Delta = \coring{C}$.
We have made the identifications $\dostensor{A}{A}{\coring{C}}
\cong \coring{C} \cong \dostensor{\coring{C}}{A}{A}$, thought that
the $A$--bimodules are assumed to be unital over $A$. We will use
Heyneman-Sweedler notation for the comultiplication, namely,
$\Delta (c) = c_{(1)} \tensor{A} c_{(2)}$ (sum is understood) for
$c \in \coring{C}$. By $\rcomod{\coring{C}}$ we denote the
category of right comodules over an $A$--coring $\coring{C}$. A
quick introduction to the basic properties of this category that
we will need may be found in \cite{Gomez:2002}. A more
comprehensive account is \cite{Brzezinski/Wisbauer:2003}.
Heyneman-Sweedler notation for comodules reads $\rho(m) = m_{[0]}
\tensor{A} m_{[1]}$ for $m$ in a right comodule $M$ with coaction
$\rho : M \rightarrow M \tensor{A} \coring{C}$.

Let $\rmd{R}$ denote the category of right modules over an
associative ring $R$. The category $\rmd{R}$ is a Grothendieck
category, being limits and colimits already computed in the
category $\mathsf{Ab}$ of abelian groups (the category $\rmd{R}$
is in fact isomorphic to the category of unital right modules over
a ring extension with unity defined in a straightforward way on
the $K$--bimodule $R \oplus K$). The tensor product functor
$\dostensor{-}{R}{R} : \rmd{R} \rightarrow \mathsf{Ab}$ is right
exact. When this functor is exact, we say that ${}_RR$ is
\emph{flat}. More generally, a left $R$--module $M$ is \emph{flat}
when the functor $\dostensor{-}{R}{M} : \rmd{R} \rightarrow
\mathsf{Ab}$ is exact.

For every right $R$--module $M$, let
\begin{equation}\label{mu}
\xymatrix{\varpi_M : \dostensor{M}{R}{R} \ar[r] & M & &
(\dostensor{m}{R}{r} \mapsto mr)}
\end{equation}
denote the ``multiplication map'', which is a homomorphism of
right $R$--modules. The product of $R$ can be understood as
$\varpi_R : \dostensor{R}{R}{R} \rightarrow R$. It is a
homomorphism of $R$--bimodules.

A right $R$--module $M$ is said to be \emph{firm} if $\varpi_M$ is
an isomorphism. When $M = R$, we will speak of a \emph{firm ring}.
According to \cite{Brzezinski/Kadison/Wisbauer:?}, this
terminology is due to D. Quillen. Firm modules are called regular
modules in \cite{Grandjean/Vitale:1998}, where their Morita theory
is developed. By $\rmod{R}$ we denote the full subcategory of
$\rmd{R}$ whose objects are the firm right $R$--modules. It is
known (\cite[Corollary
1.3]{Grandjean/Vitale:1998},\cite[Proposition 2.7]{Marin:1998})
that $\rmod{R}$ is abelian if $R$ is firm. In both cases, the
proof is indirect, namely, $\rmod{R}$ is shown to be equivalent to
a certain abelian category. In particular, due essentially to the
lack of left exactness of $\dostensor{-}{R}{R}$, kernels cannot be
in general computed in $\mathsf{Ab}$. The situation is simpler
if, in addition, ${}_RR$ is assumed to be flat: $\rmod{R}$ is then
easily shown to be a Grothendieck category, where coproducts,
kernels and cokernels are already computed in $\mathsf{Ab}$. A
generator for $\rmod{R}$ is given by the firm right $R$--module
$R$. When $R$ has a unity, then $\rmod{R}$ becomes the category of
unital right $R$--modules in the usual sense.

If $M$ is a firm right $R$--module and $d_M : M \rightarrow
\dostensor{M}{R}{R}$ is the inverse map of the multiplication map
$\varpi_R : \dostensor{M}{R}{R} \rightarrow M$, then we will use
the notation $d_M(m) = \dostensor{m^r}{R}{r} \in
\dostensor{M}{R}{R}$ (sum is understood). Of course, this element
in the tensor product is determined by the condition $m^rr = m$.
An analogous notation will be used for firm left $R$--modules.

Examples of firm rings are rings with (idempotent) local units,
such as the ones in \cite{Ab} and \cite{AM}. A comparison of
different notions of rings with local units may be found in
\cite{Vercruysse:unp2004}.

\section{Some adjoint pairs for modules}\label{pairsadj}

 Let $A$, $B$ denote rings with unity, and consider bimodules
${}_B\Sigma_A$, ${}_A\Sigma'_B$, and a homomorphism of
$A$--bimodules
\[
\mu : \Sigma' \tensor{B} \Sigma \rightarrow A.
\]
We have then a structure of $B$--ring on $S = \Sigma \tensor{A}
\Sigma'$, with the multiplication defined by the composite map
\[
\xymatrix{ \nabla : \Sigma \tensor{A} \Sigma' \tensor{B} \Sigma
\tensor{A} \Sigma' \ar^-{\mu}[r] & \Sigma \tensor{A} A \tensor{A}
\Sigma' \ar[r] & \Sigma \tensor{A} \Sigma'}.
\]
This multiplication, given explicitly on elements by the rule
\begin{equation}\label{multB}
\nabla(\fourtensor{x}{A}{\phi}{B}{y}{A}{\psi}) =
\dostensor{x\mu(\dostensor{\phi}{B}{y})}{A}{\psi} =
\dostensor{x}{A}{\mu(\dostensor{\phi}{B}{y})\psi},
\end{equation}
is easily shown to be associative. Moreover, the map $S
\rightarrow \rend{A}{\Sigma}$ that sends $\dostensor{\phi}{A}{x}$
to the endomorphism $y \mapsto x\mu(\dostensor{\phi}{B}{y})$ is a
homomorphism of $B$--rings. In particular, $\Sigma$ becomes an
$S-A$--bimodule. Analogously, the $B$--ring homomorphism $S
\rightarrow \lend{A}{\Sigma'}$ which maps $\dostensor{y}{A}{\psi}$
onto the endomorphism $\phi \mapsto \mu
(\dostensor{\phi}{B}{y})\psi$
 is a homomorphism of $B$--rings and,
thus, $\Sigma'$ becomes an $A-S$--bimodule.

A straightforward computation shows that
\begin{equation}\label{muBS}
\mu(\dostensor{\phi s}{B}{x}) = \mu(\dostensor{\phi}{B}{sx})
\end{equation}
 for
every $s \in S$. This means that we can consider the $A$--bimodule
homomorphism
\begin{equation}\label{sigmaS}
\mu : \dostensor{\Sigma'}{S}{\Sigma} \rightarrow A, \qquad (
\dostensor{\phi}{S}{x} \mapsto \mu(\dostensor{\phi}{B}{x})).
\end{equation}

Now, consider a ring homomorphism $\iota : R \rightarrow S =
\dostensor{\Sigma}{A}{\Sigma'}$, where $R$ is any ring (thus, we
are not assuming a unity). If $r \in R$, we will write $\iota (r)
= \abrir{r} \in \dostensor{\Sigma}{A}{\Sigma'}$ (sum is
understood) for its image in $S$. The left $S$--module structure
of $\Sigma$ induces a left $R$--module structure given explicitly
by the formula
\begin{equation}\label{Raccion0}
ru = e_r\mu(\dostensor{\varphi_r}{B}{u}), \qquad (r \in R, u \in
\Sigma).
\end{equation}
In this way, the $A$--bimodule map \eqref{sigmaS} induces a
homomorphism of $A$--bimodules
\begin{equation}\label{sigmaR}
\mu : \dostensor{\Sigma'}{R}{\Sigma} \rightarrow A, \qquad (
\dostensor{\phi}{R}{x} \mapsto \mu(\dostensor{\phi}{S}{x})).
\end{equation}
Therefore, the equality
\begin{equation}\label{Raccion}
ru = e_r\mu(\dostensor{\varphi_r}{R}{u}), \qquad (r \in R, u \in
\Sigma)
\end{equation}
is deduced from \eqref{Raccion0}.

\medskip

 There is an adjoint pair of
functors
\begin{equation}\label{adj1}
\xymatrix{\rmd{R} \ar@<0.5ex>^-{\dostensor{-}{R}{\Sigma}}[rr] & &
\rmod{A} \;, \ar@<0.5ex>^-{\hom{A}{\Sigma}{-}}[ll] &
\dostensor{-}{R}{\Sigma} \dashv \hom{A}{\Sigma}{-}},
\end{equation}
with unit and counit given, respectively, by
\[
\xymatrix{\nu_N : N \ar[r] &
\hom{A}{\Sigma}{\dostensor{N}{R}{\Sigma}}, & \nu_N(n)(u) =
\dostensor{n}{R}{u}}
\]
and
\[
\xymatrix{\zeta_M : \dostensor{\hom{A}{\Sigma}{M}}{R}{\Sigma}
\ar[r] & M, & \zeta(\dostensor{\phi}{R}{u}) = \phi(u)}.
\]
Now, if $R$ is a firm ring, then there is an adjoint pair
\begin{equation}\label{adj2}
\xymatrix{\rmod{R} \ar@<0.5ex>^-{J}[rr] & & \rmd{R},
\ar@<0.5ex>^{\dostensor{-}{R}{R}}[ll] & J \dashv
\dostensor{-}{R}{R}},
\end{equation}
where $J$ is the inclusion functor. The counit of this adjunction
is given by the multiplication map $\varpi_M : \dostensor{M}{R}{R}
\rightarrow M$, and the unit is given by its inverse $d_N : N
\rightarrow \dostensor{N}{R}{R}$ when $N$ is firm. As a
consequence, the functor $\dostensor{-}{R}{R}$ is left exact. We
have then that $\dostensor{-}{R}{R} : \rmd{R} \rightarrow
\mathsf{Ab}$ is exact (that is, ${}_RR$ is flat) if and only if
the inclusion functor $J$ is exact (since the forgetful functor
$\rmd{R} \rightarrow \mathsf{Ab}$ is faithful and exact).

By composing the adjunctions \eqref{adj1} and \eqref{adj2} we
obtain the adjoint pair
\begin{equation}\label{adj3}
\xymatrix{\rmod{R} \ar@<0.5ex>^-{\dostensor{-}{R}{\Sigma}}[rrr] &
& & \rmod{A},
\ar@<0.5ex>^-{\dostensor{\hom{A}{\Sigma}{-}}{R}{R}}[lll] &
\dostensor{-}{R}{\Sigma} \dashv
\dostensor{\hom{A}{\Sigma}{-}}{R}{R}}
\end{equation}
We shall need the explicit expression for its unit
\begin{equation}\label{unit3}
\xymatrix{\eta_N : N \ar[r] &
\dostensor{\hom{A}{\Sigma}{\dostensor{N}{R}{\Sigma}}}{R}{R}, &
\eta_N(n) = \dostensor{\nu_N(n^r)}{R}{r},}
\end{equation}
and its counit
\begin{equation}\label{counit3}
\xymatrix{\pi_M : \trestensor{\hom{A}{\Sigma}{M}}{R}{R}{R}{\Sigma}
\ar[r] & M, & \pi_M(\trestensor{f}{R}{r}{R}{u}) = f(ru)}.
\end{equation}
For our purposes it is interesting to recognize the former right
adjoint functor as a tensor-product functor. As in the case of
$\Sigma$, the right $S$--module structure on $\Sigma'$ induces a
structure of right $R$--module on $\Sigma'$ that satisfies
\begin{equation}\label{accionR}
\phi r = \mu(\dostensor{\phi}{R}{e_r})\varphi_r
\end{equation}

\begin{proposition}\label{homisotensor}
Let $\iota : R \rightarrow \dostensor{\Sigma}{A}{\Sigma'}$ be a
homomorphism of rings. Consider the $A-R$--bimodule $\Sigma^{\dag}
= \dostensor{\Sigma'}{R}{R}$.  If $R$ is a firm ring and $\Sigma$
is firm as a left $R$--module, then there exists a natural
isomorphism $\alpha : \dostensor{\hom{A}{\Sigma}{-}}{R}{R} \simeq
\dostensor{-}{A}{\Sigma^{\dag}}$. In particular, we have the
adjoint pair
\begin{equation}\label{adj6}
\xymatrix{\rmod{R} \ar@<0.5ex>^-{\dostensor{-}{R}{\Sigma}}[rrr] &
& & \rmod{A}, \ar@<0.5ex>^-{\dostensor{-}{A}{\Sigma^{\dag}}}[lll]
& \dostensor{-}{R}{\Sigma} \dashv \dostensor{-}{A}{\Sigma^{\dag}}}
\end{equation}
\end{proposition}
\begin{proof}
For $M \in \rmod{A}$ define
\[
\xymatrix{\alpha_M : \dostensor{\hom{A}{\Sigma}{M}}{R}{R} \ar[r] &
\dostensor{M}{A}{\Sigma^{\dag}}, & \alpha_M(\dostensor{h}{R}{r}) =
\trestensor{h(e_s)}{A}{\varphi_s}{R}{r^s}}.
\]
A straightforward computation shows that $\alpha_M$ is natural in
$M$, and that its inverse is given by
\[
\xymatrix{\alpha^{-1}_M : \dostensor{M}{R}{\Sigma^{\dag}} \ar[r] &
\dostensor{\hom{A}{\Sigma}{M}}{R}{R}, & \alpha^{-1}_M
(\trestensor{m}{A}{\phi}{R}{r}) = \dostensor{h}{R}{r},}
\]
where $h(u) = m\mu(\dostensor{\phi}{R}{u})$.
\end{proof}

It is convenient for our purposes to have an explicit expression
for the adjointness isomorphism, the unit and the counit of the
adjoint pair \eqref{adj6}. The first is given by
\begin{equation}\label{isoadj6}
\xymatrix{\hom{A}{\dostensor{N}{R}{\Sigma}}{M} \ar[r] &
\hom{R}{N}{\dostensor{M}{A}{\Sigma^{\dag}}} \\
f \ar@{|->}[r] & (n \mapsto
\trestensor{f(\dostensor{n^r}{R}{e_s})}{A}{\varphi_s}{R}{r^s}),}
\end{equation}
the counit by
\begin{equation}\label{counit6}
\xymatrix{\delta_M : \trestensor{M}{A}{\Sigma^{\dag}}{R}{\Sigma}
\ar[r] & M & (\fourtensor{m}{A}{\phi}{R}{r}{R}{x} \mapsto
m\mu(\dostensor{\phi}{R}{rx})),}
\end{equation}
and the unit by
\[
\xymatrix{\eta_M : N \ar[r] &
\trestensor{N}{R}{\Sigma}{A}{\Sigma^{\dag}} & (n \mapsto
\fourtensor{n^r}{R}{e_s}{A}{\phi_s}{R}{r^s})}
\]

\section{Comatrix corings}\label{comatrixsection}

Keep the notations of Section \ref{pairsadj}. We thus consider a
ring homomorphism $\iota : R \rightarrow S = \Sigma \tensor{A}
\Sigma'$, where the structure of ring of $\Sigma \tensor{A}
\Sigma'$ is given by \eqref{multB}. We assume that $R$ is a firm
ring, and that $A$ is a ring with unity. Recall that, for a firm
left $R$--module $M$, $d_M : M \rightarrow \dostensor{R}{R}{M}$
denotes the inverse of the multiplication map, and we use the
notation $d_{M}(m) = \dostensor{r}{R}{m^r}$.

Let $M$ be an $R-A$--bimodule which has a structure of a right
comodule over an $A$--coring $\coring{C}$. We will say that $M$ is
an $R-\coring{C}$--bicomodule if the structure map $\rho_M : M
\rightarrow \dostensor{M}{A}{\coring{C}}$ is left $R$--linear.
 The following generalizes the construction of
comatrix corings given in \cite{ElKaoutit/Gomez:2003}.

\begin{theorem}\label{comatrixconst}
If ${}_R\Sigma$ is firm as a left $R$--module, then
$\dostensor{\Sigma'}{R}{\Sigma}$ is an $A$--coring with
comultiplication
\begin{equation*}
\xymatrix{\Delta : \dostensor{\Sigma'}{R}{\Sigma}
\ar^-{\dostensor{\Sigma'}{R}{d_{\Sigma}}}[rr] & &
\trestensor{\Sigma'}{R}{R}{R}{\Sigma}
\ar^-{\trestensor{\Sigma'}{R}{\iota}{R}{\Sigma}}[rr] & &
\fourtensor{\Sigma'}{R}{\Sigma}{A}{\Sigma'}{R}{\Sigma}}
\end{equation*}
and counity
\begin{equation*}
\xymatrix{\mu : \dostensor{\Sigma'}{R}{\Sigma} \ar[r] & A}
\end{equation*}
Moreover, $\Sigma$ is an
$R-\dostensor{\Sigma'}{R}{\Sigma}$--bicomodule, and $\Sigma^{\dag}
= \dostensor{\Sigma'}{R}{R}$ is a
$\dostensor{\Sigma'}{R}{\Sigma}-R$--bicomodule with the coactions
\begin{equation*}
\xymatrix{\rho_{\Sigma} : \Sigma \ar^-{d_{\Sigma}}[r] &
\dostensor{R}{R}{\Sigma} \ar^-{\dostensor{\iota}{R}{\Sigma}}[rr] &
& \trestensor{\Sigma}{A}{\Sigma'}{R}{\Sigma}}
\end{equation*}
and
\begin{equation*}
\xymatrix{\lambda_{\Sigma^{\dag}}: \dostensor{\Sigma'}{R}{R}
\ar^-{\dostensor{\Sigma'}{R}{d_R}}[rr] & &
\trestensor{\Sigma'}{R}{R}{R}{R}
\ar^-{\trestensor{\Sigma'}{R}{\iota}{R}{R}}[rr] & &
\fourtensor{\Sigma'}{R}{\Sigma}{A}{\Sigma'}{R}{R}}
\end{equation*}
\end{theorem}
\begin{proof}
For $\dostensor{\phi}{R}{u} \in \dostensor{\Sigma'}{R}{\Sigma}$ we
have
\begin{equation}\label{comultiplication}
\Delta (\dostensor{\phi}{R}{u}) =
\trestensor{\phi}{R}{\abrir{r}}{R}{u^r}
\end{equation}
Therefore,
\[
\begin{array}{rcl}
(\trestensor{\Sigma'}{R}{\Sigma}{A}{\Delta})\Delta(\dostensor{\phi}{R}{u})
& = &
(\trestensor{\Sigma'}{R}{\Sigma}{A}{\Delta})(\trestensor{\phi}{R}{\abrir{r}}{R}{u^r})
\\
 & = & \fourtensor{\phi}{R}{\abrir{r}}{A}{\abrir{t}}{R}{u^{rt}}\\
 & = &
 (\fourtensor{\Sigma'}{R}{\iota}{R}{\iota}{R}{\Sigma})(\fourtensor{\phi}{R}{r}{R}{t}{R}{u^{rt}})\\
 & = &
 (\fourtensor{\Sigma'}{R}{\iota}{R}{\iota}{R}{\Sigma})(\fourtensor{\phi}{R}{sr^s}{R}{t}{R}{u^{rt}})\\
 & = &
 (\fourtensor{\Sigma'}{R}{\iota}{R}{\iota}{R}{\Sigma})(\fourtensor{\phi}{R}{s}{R}{r^st}{R}{u^{rt}})\\
 & = &
 (\fourtensor{\Sigma'}{R}{\iota}{R}{\iota}{R}{\Sigma})(\fourtensor{\phi}{R}{s}{R}{r^s}{R}{tu^{rt}})\\
 & = &
 (\fourtensor{\Sigma'}{R}{\iota}{R}{\iota}{R}{\Sigma})(\fourtensor{\phi}{R}{s}{R}{r^s}{R}{u^{r}}).
\end{array}
\]
On the other hand,
\[
\begin{array}{rcl}
(\trestensor{\Delta}{A}{\Sigma'}{R}{\Sigma})\Delta
(\dostensor{\phi}{R}{u}) & = &
(\trestensor{\Delta}{A}{\Sigma'}{R}{\Sigma})(\trestensor{\phi}{R}{\abrir{r}}{R}{u^r})\\
 & = &
 \fourtensor{\phi}{R}{\abrir{s}}{R}{{e_r}^s}{A}{\varphi_r}\tensor{R}
 u^r
\end{array}
\]
Thus, what we need to check to prove that $\Delta$ is
coassociative is that
\begin{equation}\label{sutil}
(\dostensor{\iota}{R}{\iota})(\dostensor{s}{R}{r^s}) =
\trestensor{\abrir{s}}{R}{{e_r}^s}{A}{\varphi_r}.
\end{equation}
To prove \eqref{sutil} let $\varpi_{\Sigma} :
\dostensor{R}{R}{\Sigma} \rightarrow \Sigma$ be the left action
map, and $\varpi_R : \dostensor{R}{R}{R} \rightarrow R$ the
multiplication map. It follows from \eqref{Raccion} that
\begin{equation*}
\nabla \circ (\iota \tensor{R} R) = \varpi_{\Sigma} \tensor{A}
\Sigma
\end{equation*}
 Since $\iota$ is an algebra map, we deduce
\begin{equation}\label{porfirm}
\iota \; \varpi_R = (\dostensor{\varpi_{\Sigma}}{A}{\Sigma'})
(\dostensor{R}{R}{\iota}).
\end{equation}
 Now, both $R$ and $\Sigma$ are firm,
which allows to rewrite \eqref{porfirm} as
\begin{equation*}
(\dostensor{d_{\Sigma}}{A}{\Sigma'}) \; \iota =
(\dostensor{R}{R}{\iota}) \; d_R
\end{equation*}
from which \eqref{sutil} follows by composing with
$\trestensor{\iota}{R}{\Sigma}{A}{\Sigma'}$ on the left. The
counitary property follows easily from \eqref{Raccion} and
\eqref{accionR}. The arguments to prove that $\rho_{\Sigma}$ and
$\lambda_{\Sigma^{\dag}}$ are comodule structure maps are similar.
To check this, it is useful to use the following expressions for
the right and left coactions:
\begin{equation}\label{sigma}
\rho_{\Sigma}(u) = \dostensor{\abrir{r}}{R}{u^r}
\end{equation}
and
\begin{equation*}
\lambda_{\Sigma^{\dag}}(\dostensor{\phi}{R}{r}) =
\trestensor{\phi}{R}{\abrir{s}}{R}{r^s}
\end{equation*}
 Obviously, they are $R$--linear.
\end{proof}

\begin{definition}\label{comatrixdefined}
We refer to the $A$--coring $\dostensor{\Sigma'}{R}{\Sigma}$ as
the \emph{comatrix coring} associated to the ring homomorphism
$\iota: R \rightarrow \Sigma \tensor{A} \Sigma'$,  the
$A$--bimodule map $\mu : \Sigma' \tensor{B} \Sigma \rightarrow A$,
and the firm module ${}_R\Sigma$.
\end{definition}

\begin{remark}
Starting from $\mu : \Sigma' \tensor{B} \Sigma \rightarrow A$ and
$\iota : R \rightarrow \Sigma \tensor{A} \Sigma'$ and assuming
$\Sigma'$ to be firm as a right $R$--module,  an  A-coring
structure is defined on $\Sigma' \tensor{R} \Sigma$ whose
comultiplication is given by
\begin{equation*}
\Delta' (\dostensor{\phi}{R}{u}) =
\trestensor{\phi^r}{R}{\abrir{r}}{R}{u},
\end{equation*}
and with counity $\mu : \Sigma' \tensor{R} \Sigma \rightarrow A$.
A straightforward version of Theorem \ref{comatrixconst} shows
that $\Sigma'$ is a left comodule, with coaction
\begin{equation}\label{sigmaprima}
\lambda_{\Sigma'}(\phi) = \dostensor{\phi^r}{R}{\abrir{r}}
\end{equation}
and ${}^{\dag}\Sigma = R \tensor{R} \Sigma$ is a right comodule,
with coaction
\begin{equation*}
\rho_{{}^{\dag}\Sigma}(\dostensor{r}{R}{u}) =
\trestensor{s}{R}{\abrir{r^s}}{R}{u},
\end{equation*}
where $r = sr^s$. Here, ${}_R\Sigma$ is not assumed to be firm.
\end{remark}

When both modules ${}_R\Sigma$ and $\Sigma'_R$ are firm, it is
equivalent to construct the comatrix coring using the
``firmeness'' of $\Sigma$ or $\Sigma'$, as the following
proposition establishes.

\begin{proposition}\label{sigmasigmaprima}
If ${}_R\Sigma$ and $\Sigma'_{R}$ are assumed to be firm, then
$\Delta = \Delta'$, and the isomorphism $d_{\Sigma'} : \Sigma'
\cong \Sigma' \tensor{R} R = \Sigma^{\dag}$ (resp. $d_{\Sigma} :
\Sigma \cong R \tensor{R} \Sigma = {}^{\dag}\Sigma$) is an
isomorphism of $\Sigma' \tensor{R} \Sigma-R$--bicomodules (resp.
of $R-\Sigma'\tensor{R}\Sigma$--bicomodules).
\end{proposition}
\begin{proof}
To prove that $\Delta = \Delta'$ it is enough to check that for
every $\dostensor{\phi}{R}{u} \in \Sigma' \tensor{R} \Sigma$ one
has
\begin{equation*}
\trestensor{\phi^r}{R}{r}{R}{u} = \trestensor{\phi}{R}{r}{R}{u^r}
\end{equation*}
But this equality follows by applying the isomorphism $\Sigma'
\tensor{R} \varpi_{\Sigma} = \varpi_{\Sigma'} \tensor{R} \Sigma$
to both members. To show that $d_{\Sigma'}$ is an isomorphism of
left $\Sigma'\tensor{R}\Sigma$--comodules we have to check the
equality
\begin{equation}\label{primadag}
\trestensor{\phi^r}{R}{\abrir{s}}{R}{r^s} =
\fourtensor{\phi^r}{R}{e_r}{A}{(\varphi_r)^s}{R}{s}
\end{equation}
Equality \eqref{primadag} holds whenever it gives a correct
equality after applying the isomorphism $\Sigma' \tensor{R} \Sigma
\tensor{A} \varpi_{\Sigma'}$. But the obtained equality then reads
\begin{equation*}
\trestensor{\phi^r}{R}{e_s}{A}{\varphi_sr^s} =
\trestensor{\phi^r}{R}{e_r}{A}{(\varphi_r)^ss},
\end{equation*}
which is correct by \eqref{Raccion}. An analogous argument shows
that $d_{\Sigma}$ is an isomorphism of $R-\Sigma' \tensor{R}
\Sigma$--bicomodules.
\end{proof}

Recall from \cite{Brzezinski/Gomez:2003} the notion of a comatrix coring context, consisting of
two rings (with unit) $A$, $B$, an $A-B$-bimodule $\Sigma'$, a $B-A$-bimodule $\Sigma$ and two bimodule maps
$$\iota:B\to\dostensor{\Sigma}{A}{\Sigma'}\quad\textrm{and}\quad \varepsilon : \dostensor{\Sigma'}{B}{\Sigma}\to A,$$
such that the following diagrams commute
\begin{equation}\label{comatrixcontext}
\xymatrix{\Sigma' \tensor{B} \Sigma \tensor{A} \Sigma'
\ar^-{\varepsilon \tensor{A} \Sigma'}[rr] & & A \tensor{A} \Sigma'
\ar[d] \\
\Sigma' \tensor{B} B \ar^-{\Sigma' \tensor{B}  \iota}[u] \ar[rr]&
& \Sigma' } \qquad \xymatrix{\Sigma \tensor{A} \Sigma' \tensor{B}
\Sigma \ar^{\Sigma \tensor{A} \varepsilon}[d] & & R \tensor{B}
\Sigma \ar^-{\iota \tensor{B} \Sigma}[ll] \ar[d] \\
\Sigma \tensor{A} A \ar[rr] & & \Sigma}
\end{equation}
It follows \cite[Theorem 2.4]{Brzezinski/Gomez:2003} that
$\dostensor{\Sigma'}{B}{\Sigma}$ is a coring, $\Sigma$ is finitely
generated and projective as a right $A$-module and $\Sigma'$ is
isomorphic to $\Sigma^*$. However, once we replace in this
construction the unital ring $B$ by a firm ring $R$, and we assume
that $\Sigma$ and $\Sigma'$ are firm as $R$-modules, we can still
construct the coring $\dostensor{\Sigma'}{R}{\Sigma}$ as before,
but there is no need anymore for $\Sigma$ to be finitely generated
and projective as right $A$-module. It is straightforward to check
that the comatrix coring as constructed from the comatrix coring
context is exactly the comatrix coring from Proposition
\ref{sigmasigmaprima}.

\begin{example}\label{fincomatrix}
If $\Sigma$ is a $B-A$--bimodule such that $\Sigma_A$ is finitely
generated and projective, and $\{ (e_i, e_i^*) \} \subseteq \Sigma
\times \Sigma^*$ is a finite dual basis for $\Sigma$, then the
canonical ring homomorphism $\iota : B \rightarrow \Sigma
\tensor{A} \Sigma^*$ given by $\iota(b) =
\dostensor{be_i}{A}{e_i^*} = \dostensor{e_i}{A}{e_i^*b}$ (sum is
understood) induces on $\rcomatrix{B}{\Sigma}$, in virtue of
Theorem \ref{comatrixconst}, a structure of $A$--coring. Here, the
$A$--bimodule map $\mu : \rcomatrix{B}{\Sigma} \rightarrow A$ is
just the evaluation map. These corings were introduced in
\cite{ElKaoutit/Gomez:2003}. We call them \emph{finite} comatrix
corings.
\end{example}

\begin{example}\label{infcomatrix}
Let $\cat{P}$ be a set of finitely generated and projective right
modules over a ring $A$, and consider $\Sigma = \bigoplus_{P \in
\cat{P}} P$.  For each $P \in \cat{P}$ let $\iota_P : P
\rightarrow \Sigma$ and $\pi_P : \Sigma \rightarrow P$ be the
canonical injection and projection, respectively. Then the
elements $u_P = \iota_P \pi_P$ with $P \in \cat{P}$ form a set of
orthogonal idempotents in $\rend{A}{\Sigma}$. Consider any
homomorphism of rings $T \rightarrow \rend{A}{\Sigma}$, and let $R
= \sum_{P,Q \in \cat{P}} u_PTu_Q$, which is a (non unital) subring
of the image of $T$ in $\rend{A}{\Sigma}$. Consider the unique
homomorphism of rings $\iota : R \rightarrow \Sigma \tensor{A}
\Sigma^*$ defined on each $u_PTu_Q$ by the composite map
\begin{equation*}
\xymatrix{u_PTu_Q \ar[r] & \hom{A}{Q}{P} \ar^-{\simeq}[r] & P
\tensor{A} Q^* \subseteq \Sigma \tensor{A} \Sigma^*,}
\end{equation*}
where the first map assigns $\pi_P f \iota_Q$ to every $f \in
u_PTu_Q$. Clearly the ring $R$ has enough idempotents, which
implies that $R$ is a firm ring. Moreover, $\Sigma$ is easily
shown to be a firm left $R$--module, with the left $R$--action
given by the restriction of its canonical structure of left
$T$--module. Therefore, the comatrix $A$--coring $\Sigma^*
\tensor{R} \Sigma$ stated in Theorem \ref{comatrixconst} makes
sense. We will see in Section \ref{comparison} that this
construction is tightly related to the \emph{infinite comatrix
corings} from \cite{ElKaoutit/Gomez:2004}.

An explicit expression for the comultiplication of
$\rcomatrix{R}{\Sigma}$ in this case is obtained as follows. For
each $P \in \cat{P}$, let $\{ e_{{\alpha}_P}, e_{{\alpha_P}}^* \}$
be a finite dual basis for $P_A$. Then
\begin{equation*}
\iota (u_P) = \sum_{\alpha_P} \iota_P(e_{\alpha_P}) \tensor{R}
e_{\alpha_P}^*\pi_P,
\end{equation*}
identity from which, in conjunction with \eqref{comultiplication},
we easily deduce, for $\phi \in \Sigma^*$ and $x \in \Sigma$, that
\begin{equation*}
\Delta(\dostensor{\phi}{R}{x}) = \sum_{P \in
\cat{F}}\fourtensor{\phi}{R}{\iota_P(e_{\alpha_P})}{A}{e_{\alpha_P}^*\pi_P}{R}{x},
\end{equation*}
where $\cat{F}$ is a finite subset of $\cat{P}$ such that $x =
\sum_{P \in \cat{F}}e_Px$. The $\rcomatrix{R}{\Sigma}$--comodule
structures of $\Sigma$ and $\Sigma^{\dag} = \Sigma \tensor{R} R$
are then given by
\begin{equation*}
\rho_{\Sigma}(x)  = \sum_{P \in
\cat{F}}\trestensor{\iota_P(e_{\alpha_P})}{A}{e_{\alpha_P}^*\pi_P}{R}{x},
\end{equation*}
and
\begin{equation*}
\lambda_{\Sigma^{\dag}} (\dostensor{\phi}{R}{r}) = \sum_{Q \in
\cat{G}}
\fourtensor{\phi}{R}{\iota_Q(e_{\beta_Q})}{A}{e_{\alpha_Q}^*\pi_Q}{R}{r},
\end{equation*}
respectively. In the second case, $\cat{G}$ is a finite subset of
$\cat{P}$ such that $r = \sum_{Q \in \cat{P}} e_Qr$.
\end{example}

\begin{example}\label{locprojcom}
Let $\Sigma$ be a $B-A$-bimodule that is locally projective in the
following sense (see \cite{Vercruysse:unp2004}): for every element
$u\in\Sigma$ (or equivalently, for every finite set in $\Sigma$),
we can find an element (called dual basis)
$e=\dostensor{e_i}{A}{f_i}\in\dostensor{\Sigma}{A}{\Sigma^*}$ such
that $be=eb$ and $e\cdot u=u$. In case where $B=\mathbb{Z}$, this
is equivalent with local projectivity in the sense of
Zimmermann-Huisgen \cite{Zimmermann:1976}. Take now a set of
generators $U$ for $\Sigma_A$, and let $E$ be a set consisting of
local dual bases for each of the elements of $U$. Put $R$ the
subring of $\dostensor{\Sigma}{A}{\Sigma^*}$ generated by the
elements of $E$. By construction, $R$ acts with local units on
$\Sigma$ and thus $\Sigma$ is a firm left $R$-module.  Clearly,
$R$ is a ring with left local units and, therefore, it is a firm
ring (see also \cite[Theorem 3.4]{Vercruysse:unp2004}). It follows
now from our general theory (Theorem \ref{comatrixconst}) that
$\dostensor{\Sigma^*}{R}{\Sigma}$ is an $A$-coring. The counit and
comultiplication are given explicitly by the following formulas
\begin{equation*}
\varepsilon(\dostensor{\varphi}{R}{u})=\varphi(u).
\end{equation*}
and
\begin{equation*}
\Delta(\dostensor{\varphi}{R}{u})=\fourtensor{\varphi}{R}{e_i}{A}{f_i}{R}{u}=\trestensor{\varphi}{R}{e}{R}{u},
\end{equation*}
where $e=\dostensor{e_i}{A}{f_i}$ is a dual basis for $u$. That
this comultiplication is well-defined can be also checked directly
as follows. If $e'$ is another dual basis for $u$, let then $e''$
be a local unit for both $e$ and $e'$, i.e. $e''e=e$ and
$e''e'=e'$. We compute
\begin{equation*}
\begin{array}{rcl}
\trestensor{\varphi}{R}{e}{R}{u}&=&\trestensor{\varphi}{R}{e''e}{R}{u}=\trestensor{\varphi}{R}{e''}{R}{eu}\\
&=&\trestensor{\varphi}{R}{e''}{R}{u}=\trestensor{\varphi}{R}{e''}{R}{e'u}\\
&=&\trestensor{\varphi}{R}{e''e'}{R}{u}=\trestensor{\varphi}{R}{e'}{R}{u}.
\end{array}
\end{equation*}
\end{example}

The comatrix coring built in Theorem \ref{comatrixconst} admits
(up to isomorphism) two more natural constructions. The first of
them uses $\Sigma^{\dag} \tensor{R} \Sigma$ as underlying
$A$--bimodule. The idea is to use as initial data the homomorphism
of $A$--bimodules $\varepsilon^{\dag} : \Sigma^{\dag} \tensor{R}
\Sigma \rightarrow A$ defined by
\begin{equation*}
\varepsilon^{\dag}(\trestensor{\phi}{R}{r}{R}{x}) = \mu
(\dostensor{\phi}{R}{rx}) = \mu(\dostensor{\phi}{B}{rx})
\end{equation*}
and the map $\iota' : R \rightarrow \Sigma \tensor{A}
\Sigma^{\dag}$ given by
\begin{equation*}
\iota'(r) = (\iota \tensor{R} R)d_R(r) = e_s \tensor{A} \varphi_s
\tensor{R} r^s
\end{equation*}

\begin{lemma}
If we endow $\Sigma \tensor{A} \Sigma^{\dag}$ with the ring
structure induced by $\varepsilon^{\dag}$, then $\iota'$ is a
homomorphism of rings.
\end{lemma}
\begin{proof}
Given $r,t \in R$, we compute
\begin{equation*}
\begin{array}{lcl}
\iota'(r)\iota'(t) & = &
(\dostensor{\abrir{s}}{R}{r^s})(\dostensor{\abrir{u}}{R}{t^u}) \\
 & = &
 \trestensor{e_s\mu(\dostensor{\varphi_s}{R}{r^se_u})}{A}{\varphi_u}{R}{t^u}
 \\
 & = & \trestensor{sr^se_u}{A}{\varphi_u}{R}{t^u} \\
 & = & \trestensor{re_u}{A}{\varphi_u}{R}{t^u} \\
 & = &
 \trestensor{e_r\mu(\dostensor{\varphi_r}{R}{e_u})}{A}{\varphi_u}{R}{t^u}\\
 & = & \dostensor{(\abrir{r})(\abrir{u})}{R}{t^u}\\
 & = & \trestensor{e_{ru}}{A}{\varphi_{ru}}{R}{t^u} \\
 & = & \trestensor{e_v}{A}{\varphi_v}{R}{(rt)^v} \\
 & = & \iota'(rt)
\end{array}
\end{equation*}
\end{proof}

In view of this lemma, we can construct the comatrix coring
associated to the firm module ${}_R\Sigma$ (or, in virtue of
Proposition \ref{sigmasigmaprima} the firm module $\Sigma^{\dag}_R
= \Sigma' \tensor{R} R$). This comatrix $A$--coring is the
$A$--bimodule $\Sigma^{\dag} \tensor{R} \Sigma$ with the
comultiplication
\begin{equation}\label{deltadag}
\Delta^{\dag}(\trestensor{\phi}{R}{r}{R}{x}) =
\fourtensor{\phi}{R}{r}{R}{\abrir{t}}{R}{s^t}\tensor{R}{x^s}
\end{equation}
and counity
\begin{equation}
\varepsilon^{\dag}(\trestensor{\phi}{R}{r}{R}{x}) = \mu
(\dostensor{\phi}{R}{rx})
\end{equation}

By Proposition \ref{sigmasigmaprima},
\begin{equation}\label{Deltadag2}
\Delta^{\dag}(\trestensor{\phi}{R}{r}{R}{x}) =
\fourtensor{\phi}{R}{s}{R}{\abrir{t}}{R}\dostensor{(r^s)^t}{R}{x}
\end{equation}
The structure of right $\Sigma^{\dag} \tensor{R} \Sigma$--comodule
on $\Sigma$ is given by (see \eqref{sigma})
\begin{equation}\label{rodag}
\rho^{\dag}_{\Sigma} (x) = \trestensor{\abrir{s}}{R}{r^s}{R}{x^r}
\end{equation}
 Now, being $\Sigma^{\dag} = \Sigma'
\tensor{R} R$ firm as a right $R$--module, it has a structure of
left $\Sigma^{\dag} \tensor{R} \Sigma$--comodule with coaction
(see \eqref{sigmaprima})
\begin{equation}\label{lambdadag}
\lambda^{\dag}_{\Sigma^{\dag}}(\dostensor{\phi}{R}{r}) =
\fourtensor{\phi}{R}{s}{R}{\abrir{t}}{R}{(r^s)^t}
\end{equation}

\begin{proposition}\label{comparar}
The map $f : \Sigma^{\dag} \tensor{R} \Sigma \rightarrow \Sigma'
\tensor{R} \Sigma$ given by $f(\trestensor{\phi}{R}{r}{R}{x}) =
\dostensor{\phi}{R}{rx}$ is an isomorphism of $A$--corings.
Furthermore, the structure of right (resp. left) comodule over
$\Sigma^{\dag}\tensor{R}\Sigma$ of $\Sigma$ (resp.
$\Sigma^{\dag}$) given in \eqref{rodag} (resp. \eqref{lambdadag})
correspond under this isomorphism of categories to the $\Sigma'
\tensor{R} \Sigma$--comodule structures established in Theorem
\ref{comatrixconst}.
\end{proposition}
\begin{proof}
First, observe that $f$ is nothing but the inverse of the
isomorphism $\dostensor{\Sigma'}{R}{d_{\Sigma}}$. The
comultiplication of $\Sigma' \tensor{R} \Sigma$ was given in
\eqref{comultiplication}. To prove that $f$ is a homomorphism of
$A$--corings we need to check that
\begin{equation}\label{fcoring}
\fourtensor{\phi}{R}{re_t}{A}{\varphi_t}{R}{s^tx^s} =
\fourtensor{\phi r}{R}{e_s}{A}{\varphi_s}{R}{x^s}
\end{equation}
for every $\trestensor{\phi}{R}{r}{R}{x} \in
\trestensor{\Sigma'}{R}{R}{R}{\Sigma}$. The equality
\eqref{fcoring} is a direct consequence of
\begin{equation}\label{colapsa}
\dostensor{e_t}{A}{\varphi_ts^t} = \dostensor{e_s}{A}{\varphi_s},
\end{equation}
and this last is deduced from the fact that the structure of right
$R$--module of $\Sigma'_R$ is given by \eqref{accionR}. Finally,
let us check the claims concerning the comodule structures. Apply
$\Sigma \tensor{A} f$ to \eqref{rodag} to obtain, in view of
\eqref{colapsa}, the induced $\Sigma' \tensor{R} \Sigma$-coaction
\begin{equation*}
(\Sigma \tensor{A} f)\rho^{\dag}_{\Sigma} (x) =
\dostensor{\abrir{r}}{R}{x^r},
\end{equation*}
which proves that $\rho_{\Sigma} = (\Sigma \tensor{A}
f)\rho^{\dag}_{\Sigma}$. As for $\Sigma^{\dag}$ concerns, the
$\Sigma' \tensor{R} \Sigma$--comodule structure map induced by the
coring isomorphism $f$ is then given by
\begin{equation*}
\begin{array}{lcl}
(f \tensor{A}
\Sigma^{\dag})(\fourtensor{\phi}{R}{s}{R}{\abrir{t}}{R}{(r^s)^t})
& = & \fourtensor{\phi}{R}{se_t}{A}{\varphi_t}{R}{(r^s)^t} \\
& = &
\fourtensor{\phi}{R}{e_s\mu(\dostensor{\varphi_s}{R}{e_t})}{A}{\varphi_t}{R}{(r^s)^t}
\\
& = &
\fourtensor{\phi}{R}{e_s}{A}{\mu(\dostensor{\varphi_s}{R}{e_t})\varphi_t}{R}{(r^s)^t}
\\
& = & \fourtensor{\phi}{R}{e_s}{A}{\varphi_st}{R}{(r^s)^t}\\
& = & \trestensor{\phi}{R}{\abrir{s}}{R}{r^s}
\end{array}
\end{equation*}
In the above computation we made use of \eqref{Raccion} and
\eqref{accionR}. The resulting map is precisely the coaction
$\lambda_{\Sigma^{\dag}}$ given by Theorem \ref{comatrixconst}.
\end{proof}

Finite comatrix corings, as defined in Example \ref{fincomatrix},
and the corings defined in Example \ref{infcomatrix}, are built on
the canonical evaluation map $\mu : \rcomatrix{B}{\Sigma}
\rightarrow A$. Our next aim is to prove that this is always the
case in some precise sense, namely a third construction of the
comatrix coring built on the $A$--bimodule $\Sigma^* \tensor{R}
\Sigma$ (see Proposition \ref{coringiso}).

\medskip

For every $\phi \in \Sigma'$ let $\overline{\phi} : \Sigma
\rightarrow A$ be defined by $\overline{\phi}(u) =
\mu(\dostensor{\phi}{R}{u})$. This gives a homomorphism
$A-R$--bimodules $\overline{(-)} : \Sigma' \rightarrow \Sigma^*$
(which is also of $A-B$--bimodules). Now, let
$\epsilon:\dostensor{\Sigma^*}{B}{\Sigma} \rightarrow A$ be the
canonical evaluation map, and consider on
$\dostensor{\Sigma}{A}{\Sigma^*}$ the corresponding ring structure
with multiplication defined by
$(\dostensor{x}{A}{\phi})(\dostensor{y}{A}{\nu}) =
\dostensor{x\phi(y)}{A}{\nu}$. Then the map
\begin{equation}\label{ringmorphism}
\xymatrix{\dostensor{\Sigma}{A}{\overline{(-)}}:
\dostensor{\Sigma}{A}{\Sigma'} \ar[r] &
\dostensor{\Sigma}{A}{\Sigma^*}}
\end{equation}
is a ring homomorphism. Moreover, the left
$\dostensor{\Sigma}{A}{\Sigma^*}$--module structure of $\Sigma$ is
compatible with this ring morphism in the sense that it induces,
by restriction of scalars, the old
$\dostensor{\Sigma}{A}{\Sigma'}$--module structure. If $\iota : R
\rightarrow \dostensor{\Sigma}{A}{\Sigma'}$ is a homomorphism of
rings, then its composition with \eqref{ringmorphism} gives a
homomorphism of rings from $R$ to
$\dostensor{\Sigma}{A}{\Sigma^*}$ which induces on $\Sigma$ the
same left $R$--module structure than $\iota$. Thus, no ambiguity
arises if we write $\iota : R \rightarrow
\dostensor{\Sigma}{A}{\Sigma^*}$. Accordingly to Theorem
\ref{comatrixconst} we can consider the $A$--corings
$\dostensor{\Sigma'}{R}{\Sigma}$ and
$\dostensor{\Sigma^*}{A}{\Sigma}$.

\begin{proposition}\label{coringiso}
The map
\begin{equation}
\xymatrix{\dostensor{\overline{(-)}}{R}{\Sigma} :
\dostensor{\Sigma'}{R}{\Sigma} \ar[r] &
\dostensor{\Sigma^*}{R}{\Sigma}}
\end{equation}
is an isomorphism of $A$--corings.
\end{proposition}
\begin{proof}
Clearly, the map is a homomorphism of $A$--bimodules. Moreover, an
easy computation shows that
$\dostensor{\overline{(-)}}{R}{\Sigma}$ is a homomorphism of
$A$--corings. We need thus to prove that it is bijective. By
Proposition \ref{homisotensor} we have an isomorphism of
$A-R$--bimodules $\alpha_A : \dostensor{\Sigma^*}{R}{R}
\rightarrow \dostensor{A}{A}{\Sigma^{\dag}}$. We obtain an
isomorphism of $A$--bimodules $\widetilde{\alpha}_A :
\dostensor{\Sigma^*}{R}{\Sigma} \cong
\dostensor{\Sigma'}{R}{\Sigma}$ as the composition
\[
\xymatrix{ \dostensor{\Sigma^*}{R}{\Sigma}
\ar^-{\dostensor{\Sigma^*}{R}{d_{\Sigma}}}[rr] & &
\trestensor{\Sigma^*}{R}{R}{R}{\Sigma}
\ar^-{\dostensor{\alpha_A}{R}{\Sigma}}[rr] & &
\trestensor{\Sigma'}{R}{R}{R}{\Sigma}
\ar^-{\dostensor{\Sigma'}{R}{\varpi_{\Sigma}}}[rr] & &
\dostensor{\Sigma'}{R}{\Sigma}}.
\]
This isomorphism acts on elements by
\[
\widetilde{\alpha}_A (\dostensor{\phi}{R}{u}) =
\dostensor{\phi(e_s)\varphi_s}{R}{r^su^r} \qquad
(\dostensor{\phi}{R}{u} \in \dostensor{\Sigma^*}{R}{\Sigma}).
\]
Now, the computation
\begin{multline*}
(\dostensor{\overline{(-)}}{R}{\Sigma}) \widetilde{\alpha}_A
(\dostensor{\phi}{R}{u}) =
\dostensor{\phi(e_s)\overline{\varphi_s}}{R}{r^su^r} = \\
\dostensor{\phi(e_s)\overline{\varphi_s}r^s}{R}{u^r} =
\dostensor{\phi r}{R}{u^r} = \dostensor{\phi}{R}{ru^r} =
\dostensor{\phi}{R}{u}
\end{multline*}
proves that $\dostensor{\overline{(-)}}{R}{\Sigma}$ is a left
inverse to the bijection $\widetilde{\alpha}_A$ and, henceforth,
$\dostensor{\overline{(-)}}{R}{\Sigma} =
\widetilde{\alpha}_A^{-1}$.
\end{proof}

\begin{remark}
The isomorphism of $A-R$--bimodules $\alpha_A :
\dostensor{\Sigma^*}{R}{R} \cong \trestensor{A}{A}{\Sigma'}{R}{R}
\cong \dostensor{\Sigma'}{R}{R}$ deduced from Proposition
\ref{homisotensor} induces a natural isomorphism between the
functors $\trestensor{-}{A}{\Sigma^*}{R}{R},
\trestensor{-}{A}{\Sigma'}{R}{R} : \rmod{A} \rightarrow \rmod{R}$.
This, together with Proposition \ref{coringiso} allows to replace
$\Sigma'$ by $\Sigma^*$ and $\mu$ by $\epsilon$ in our theory.
\end{remark}

\section{Galois Comodules}\label{Galoiscomodules}

Let $(\coring{C},\Delta_{\coring{C}},\epsilon_{\coring{C}})$ be a
coring over $A$ and $\Sigma$ be a right $\coring{C}$--comodule.
Let $\iota : R \rightarrow \dostensor{\Sigma}{A}{\Sigma^*}$ be a
ring homomorphism, for $R$ a firm ring, such that $\Sigma$ is an
$R -\coring{C}$--bicomodule and ${}_R\Sigma$ is firm. Our first
objective is to prove that the adjoint pair \eqref{adj1} induces
an adjoint pair
\begin{equation}\label{adj4a}
\xymatrix{\rmd{R} \ar@<0.5ex>^-{\dostensor{-}{R}{\Sigma}}[rrr] & &
& \rcomod{\coring{C}},
\ar@<0.5ex>^-{\hom{\coring{C}}{\Sigma}{-}}[lll] &
\dostensor{-}{R}{\Sigma} \dashv \hom{\coring{C}}{\Sigma}{-}}.
\end{equation}
\begin{lemma}\label{adjtilda}
Let $\tau : \hom{A}{\dostensor{-}{R}{\Sigma}}{-} \rightarrow
\hom{R}{-}{\hom{A}{\Sigma}{-}}$ denote the adjunction isomorphism
for the adjoint pair \eqref{adj1}. Then for every $N \in \rmd{R}$
and every $M \in  \rcomod{\coring{C}}$ we have a commutative
square
\[
\xymatrix{\hom{A}{\dostensor{N}{R}{\Sigma}}{M}
\ar^-{\tau_{N,M}}[rr] & &
\hom{R}{N}{\hom{A}{\Sigma}{M}} \\
\hom{\coring{C}}{\dostensor{N}{R}{\Sigma}}{M}
\ar^-{\tau_{N,M}}[rr] \ar[u]& &
\hom{R}{N}{\hom{\coring{C}}{\Sigma}{M}}\ar[u],}
\]
where the vertical arrows are inclusions. In particular,
$\dostensor{-}{R}{\Sigma}$ is left adjoint to
$\hom{\coring{C}}{\Sigma}{-}$.
\end{lemma}
\begin{proof}
Consider the diagram

\begin{equation}\label{equalizers1}
\xymatrix{\hom{\coring{C}}{\dostensor{N}{R}{\Sigma}}{M}  \ar[d] &
&
\hom{R}{N}{\hom{\coring{C}}{\Sigma}{M}} \ar[d] \\
\hom{A}{\dostensor{N}{R}{\Sigma}}{M} \ar[rr]^-{\tau_{N,M}}
\ar@<1ex>[d]^{\delta} \ar@<-1ex>[d]_{\beta} & &
\hom{R}{N}{\hom{A}{\Sigma}{M}} \ar@<1ex>[d]^{\delta'}
\ar@<-1ex>[d]_{\beta'} \\
\hom{A}{\dostensor{N}{R}{\Sigma}}{\dostensor{M}{A}{\coring{C}}}
\ar[rr]^-{\tau_{N,\dostensor{M}{A}{\coring{C}}}} & &
\hom{R}{N}{\hom{A}{\Sigma}{\dostensor{M}{A}{\coring{C}}}}}
\end{equation}
where the vertical sides represent equalizer diagrams. Here,
$\beta = \hom{A}{\dostensor{N}{R}{\Sigma}}{\rho_M}$ and $\beta' =
\hom{R}{N}{\hom{A}{\Sigma}{\rho_M}}$, whence the naturality of
$\tau$ gives that $\beta' \circ \tau_{N,M} =
\tau_{N,\dostensor{M}{A}{\coring{C}}} \circ \beta$. Obviously, the
kernel of $\beta - \delta$ is
$\hom{\coring{C}}{\dostensor{N}{R}{\Sigma}}{M}$. On the other
hand, $\delta (f) =
(\dostensor{f}{A}{\coring{C}})(\dostensor{N}{A}{\rho_{\Sigma}})$
for $f \in \hom{A}{\dostensor{N}{R}{\Sigma}}{M}$, and
$\delta'(g)(n) = (\dostensor{g(n)}{A}{\coring{C}})\rho_{\Sigma}$
for $g \in \hom{R}{N}{\hom{A}{\Sigma}{M}}$ and $n \in N$. Now,
given $f \in \hom{A}{\dostensor{N}{R}{\Sigma}}{M}$ and $n \in N$,
we get
\begin{equation}\label{cuenta1}
\begin{array}{rcl}
\tau_{N,\dostensor{M}{A}{\coring{C}}}(\delta(f))(n) & = &
\hom{A}{\Sigma}{\delta(f)}\nu_N(n) \\
 & = & \delta(f)\nu_N(n) \\
 & = & (\dostensor{f}{A}{\coring{C}})(\dostensor{N}{A}{\rho_{\Sigma}})\nu_N(n)
\end{array}
\end{equation}
If we apply $\delta' \circ \tau_{N,M}$ then we obtain
\begin{equation}\label{cuenta2}
\begin{array}{rcl}
\delta' (\tau_{N,M}(f))(n) & = & \delta'
(\hom{A}{\Sigma}{f}\nu_N)(n) \\
 & = &
 (\dostensor{\hom{A}{\Sigma}{f}\nu_N(n)}{A}{\coring{C}})\rho_{\Sigma}
 \\
 & = & (\dostensor{f\nu_N(n)}{A}{\coring{C}})\rho_{\Sigma}
\end{array}
\end{equation}
The values of the maps \eqref{cuenta1} and \eqref{cuenta2} at $u
\in \Sigma$ are
\begin{equation*}
\begin{array}{rcl}
(\dostensor{f}{A}{\coring{C}})(\dostensor{N}{A}{\rho_{\Sigma}})\nu_N(n)(u)
 & = &
 (\dostensor{f}{A}{\coring{C}})(\dostensor{N}{A}{\rho_{\Sigma}})(\dostensor{n}{A}{u})\\
  & = &
  (\dostensor{f}{A}{\coring{C}})(\trestensor{n}{A}{u_{[0]}}{A}{u_{[1]}})
  \\
  & = & \dostensor{f(\dostensor{n}{A}{u_{[0]}})}{A}{u_{[1]}},
\end{array}
\end{equation*}
and
\begin{equation*}
\begin{array}{rcl}
(\dostensor{f\nu_N(n)}{A}{\coring{C}})\rho_{\Sigma}(u) & = &
(\dostensor{f\nu_N(n)}{A}{\coring{C}})(\dostensor{u_{[0]}}{A}{u_{[1]}}) \\
 & = & \dostensor{f\nu_N(n)(u_{[0]})}{A}{u_{[1]}} \\
 & = & \dostensor{f(\dostensor{n}{A}{u_{[0]}})}{A}{u_{[1]}},
\end{array}
\end{equation*}
respectively, and, hence, they are equal. This concludes a proof
of the equality $\tau_{N,\dostensor{M}{A}{\coring{C}}} \circ
\delta = \delta' \circ \tau_{N,M}$. It is easy to check that the
maps $\beta'$ and $\delta'$ come out by applying the functor
$\hom{R}{N}{-}$ to the equalizer
\begin{equation}\label{equalizer3}
\xymatrix{\hom{\coring{C}}{\Sigma}{M} \ar[r] & \hom{A}{\Sigma}{M}
\ar@<0.5ex>[r]^-{\beta''} \ar@<-0.5ex>[r]_-{\delta''} &
\hom{A}{\Sigma}{\dostensor{M}{A}{\coring{C}}}},
\end{equation}
where $\beta'' = \hom{A}{\Sigma}{\rho_M}$ and $\delta''(-) =
(\dostensor{-}{A}{\coring{C}})\rho_{\Sigma}$.  Thus the equalizer
of $\beta'$ and $\delta'$ is
$\hom{R}{N}{\hom{\coring{C}}{\Sigma}{M}}$. This proves that
$\tau_{N,M}$ induces an isomorphism
\[
\tau_{N,M} : \hom{\coring{C}}{\dostensor{N}{R}{\Sigma}}{M} \simeq
\hom{R}{N}{\hom{\coring{C}}{\Sigma}{M}}.
\]
\end{proof}

Now let us, with the help of Lemma \ref{adjtilda}, prove that the
adjoint pair \eqref{adj3} induces an adjoint pair

\begin{equation}\label{adj4}
\xymatrix{\rmod{R} \ar@<0.5ex>^-{\dostensor{-}{R}{\Sigma}}[rrr] &
& & \rcomod{\coring{C}},
\ar@<0.5ex>^-{\dostensor{\hom{\coring{C}}{\Sigma}{-}}{R}{R}}[lll]
& \dostensor{-}{R}{\Sigma} \dashv
\dostensor{\hom{\coring{C}}{\Sigma}{-}}{R}{R}}.
\end{equation}
\begin{lemma}
Let $\gamma : \hom{A}{\dostensor{-}{R}{\Sigma}}{-} \rightarrow
\hom{R}{-}{\dostensor{\hom{A}{\Sigma}{-}}{R}{R}}$ denote the
adjunction isomorphism for the adjoint pair \eqref{adj3}. Then for
every $N \in \rmod{R}$ and every $M \in  \rcomod{\coring{C}}$ we
have a commutative square
\[
\xymatrix{\hom{A}{\dostensor{N}{R}{\Sigma}}{M}
\ar^-{\gamma_{N,M}}[rr] & &
\hom{R}{N}{\dostensor{\hom{A}{\Sigma}{M}}{R}{R}} \\
\hom{\coring{C}}{\dostensor{N}{R}{\Sigma}}{M}
\ar^-{\gamma_{N,M}}[rr] \ar[u]& &
\hom{R}{N}{\dostensor{\hom{\coring{C}}{\Sigma}{M}}{R}{R}} \ar[u],}
\]
where the vertical arrows are canonical monomorphisms (recall that
$\dostensor{-}{R}{R} : \rmd{R} \rightarrow \rmod{R}$ is left
exact). In particular, $\dostensor{-}{R}{\Sigma}$ is left adjoint
to $\dostensor{\hom{\coring{C}}{\Sigma}{-}}{R}{R}$.
\end{lemma}
\begin{proof}
 Consider the functorial isomorphsim
\begin{equation*}
\varsigma : \hom{R}{-}{\hom{A}{\Sigma}{-}} \simeq
\hom{R}{-}{\dostensor{\hom{A}{\Sigma}{-}}{R}{R}},
\end{equation*}
 defined from the adjunction isomorphism $\theta_{-,-} : \hom{R}{-}{-} \simeq \hom{R}{-}{\dostensor{-}{R}{R}}$ corresponding to
the adjoint pair \eqref{adj2} as $\varsigma =
\theta_{-,\hom{A}{\Sigma}{-}}$. Since $\varsigma$ is natural, we
get from the equalizer \eqref{equalizer3} a commutative diagram of
equalizers
\begin{equation*}
\xymatrix{\hom{R}{N}{\hom{\coring{C}}{\Sigma}{M}} \ar[d] & &
\hom{R}{N}{\dostensor{\hom{\coring{C}}{\Sigma}{M}}{R}{R}} \ar[d] \\
\hom{R}{N}{\hom{A}{\Sigma}{M}} \ar@<1ex>[d] \ar@<-1ex>[d]
\ar[rr]^-{\varsigma_{N,M}}& &
\hom{R}{N}{\dostensor{\hom{A}{\Sigma}{M}}{R}{R}} \ar@<1ex>[d] \ar@<-1ex>[d] \\
\hom{R}{N}{\hom{A}{\Sigma}{\dostensor{M}{A}{\coring{C}}}}
\ar[rr]^-{\varsigma_{N,\dostensor{M}{A}{\coring{C}}}}& &
\hom{R}{N}{\dostensor{\hom{A}{\Sigma}{\dostensor{M}{A}{\coring{C}}}}{R}{R}}}
\end{equation*}
which induces an isomorphism
\[
\varsigma_{N,M} : \hom{R}{N}{\hom{\coring{C}}{\Sigma}{M}} \simeq
\hom{R}{N}{\dostensor{\hom{\coring{C}}{\Sigma}{M}}{R}{R}}.
\]
Finally, we have that $\gamma = \varsigma \circ \tau$, which
finishes the proof.
\end{proof}

Let
\begin{equation}\label{counit4}
\xymatrix{\pi :
\trestensor{\hom{\coring{C}}{\Sigma}{-}}{R}{R}{R}{\Sigma} \ar[r] &
id_{\rcomod{\coring{C}}}},
\end{equation}
be the counit of the adjunction \eqref{adj4}, which is, when
evaluated at $M$, nothing but the restriction of the map displayed
in \eqref{counit3}.

When the right $\coring{C}$--comodule $\Sigma$ is finitely
generated and projective as a right $A$--module, the finite
comatrix coring $\rcomatrix{T}{\Sigma}$ is related to the coring
$\coring{C}$ with a canonical homomorphism of $A$--corings
\cite[Proposition 2.7]{ElKaoutit/Gomez:2003}, where $T =
\rend{\coring{C}}{\Sigma}$. Expressed in Sweedler's
sigma-notation, that canonical map reads
$\mathsf{can}(\dostensor{\phi}{T}{x}) = \phi(x_{[0]})x_{[1]}$. The
following proposition shows that this canonical map stills to be a
homomorphism of corings in our present more general setting.

\begin{proposition}\label{can1}
The map
\[
\xymatrix{\mathsf{can} : \dostensor{\Sigma^*}{R}{\Sigma} \ar[r] &
\coring{C}, & \dostensor{\phi}{R}{x} \mapsto \phi(x_{[0]})x_{[1]}}
\]
is a homomorphism of $A$--corings.
\end{proposition}
\begin{proof}
Given $\dostensor{\phi}{A}{x} \in \rcomatrix{R}{\Sigma}$, we have
\begin{equation*}
\begin{array}{rcl}
\phi(x_{[0]})x_{[1]} & = &
(\dostensor{\phi}{A}{\coring{C}})(\dostensor{x_{[0]}}{A}{x_{[1]}}) \\
& = & (\dostensor{\phi}{A}{\coring{C}})\rho_M(x) \\
& = & (\dostensor{\phi}{A}{\coring{C}})\rho_M(rx^r) \\
& = & (\dostensor{\phi}{A}{\coring{C}})r\rho_M(x^r) \\
& = &
(\dostensor{\phi}{A}{\coring{C}})(\dostensor{r{x^r}_{[0]}}{A}{{x^r}_{[1]}})
\\
& = & \phi(r{x^r}_{[0]}){x^r}_{[1]}.
\end{array}
\end{equation*}
Therefore,
\begin{equation}\label{cancompli}
\mathsf{can}(\dostensor{\phi}{A}{x}) =
\phi(r{x^r}_{[0]}){x^r}_{[1]}.
\end{equation}
Now, apply \eqref{cancompli} to compute
\begin{equation}\label{cancompli1}
\begin{array}{rcl}
(\Delta_{\coring{C}} \circ \mathsf{can})(\dostensor{\phi}{A}{x}) &
= &
\Delta_{\coring{C}}(\phi(r{x^r}_{[0]}){x^r}_{[1]}) \\
& = & \phi(r{x^r}_{[0]})\Delta_{\coring{C}}({x^r}_{[1]}) \\
& = &
(\trestensor{\phi}{A}{\coring{C}}{A}{\coring{C}})(\dostensor{r{x^r}_{[0]}}{A}{\Delta_{\coring{C}}({x^r}_{[1]})}).
\end{array}
\end{equation}
Let us compute the other pertinent composition:
\begin{equation}\label{cancompli2}
\begin{array}{rcl}
(\dostensor{\mathsf{can}}{A}{\mathsf{can}})\Delta_{\dostensor{\Sigma^{*}}{R}{\Sigma}}(\dostensor{\phi}{A}{x})
& = &
(\dostensor{\mathsf{can}}{A}{\mathsf{can}})(\trestensor{\phi}{A}{\abrir{r}}{R}{x^r})
\\
& = &
\dostensor{\phi({e_r}_{[0]}){e_r}_{[1]}}{A}{\varphi_r({x^r}_{[0]}){x^r}_{[1]}} \\
& = &
\dostensor{\phi({e_r}_{[0]}){e_r}_{[1]}\varphi_r({x^r}_{[0]})}{A}{{x^r}_{[1]}} \\
& = &
(\trestensor{\phi}{A}{\coring{C}}{A}{\coring{C}})(\trestensor{{e_r}_{[0]}}{A}{{e_r}_{[1]}\varphi_r({x^r}_{[0]})}{A}{{x^r}_{[1]}})
\\
& = &
(\trestensor{\phi}{A}{\coring{C}}{A}{\coring{C}})(\dostensor{\rho_{\Sigma}(e_r\varphi_r({x^r}_{[0]}))}{A}{{x^r}_{[0]}})
\\
& = &
(\trestensor{\phi}{A}{\coring{C}}{A}{\coring{C}})(\dostensor{\rho_{\Sigma}(r{x^r}_{[0]})}{A}{{x^r}_{[1]}}).
\end{array}
\end{equation}
The last members of \eqref{cancompli1} and \eqref{cancompli2} are
equal because $\Sigma$ is a right $\coring{C}$--comodule. We have
thus proved that $\Delta_{\coring{C}} \circ \mathsf{can} =
(\dostensor{\mathsf{can}}{A}{\mathsf{can}})\circ
\Delta_{\dostensor{\Sigma^{*}}{R}{\Sigma}}$. A straightforward
argument shows that $\epsilon_{\coring{C}} \circ \mathsf{can}$ is
the evaluation map, and thus the proof is complete.
\end{proof}

The homomorphism of $A$--corings $\mathsf{can}$ defined in
Proposition \ref{can1} induces two functors $\mathsf{CAN_r} :
\rcomod{\dostensor{\Sigma^*}{R}{\Sigma}} \rightarrow
\rcomod{\coring{C}}$ and $\mathsf{CAN}_l :
\lcomod{\dostensor{\Sigma^*}{R}{\Sigma}} \rightarrow
\lcomod{\coring{C}}$ which are examples of induction functors in
the sense of \cite{Gomez:2002}. In this case, the functors act as
the identity on morphisms and on the underlying structures of
$A$--modules. The image under $\mathsf{CAN}_l$ of the right
$\dostensor{\Sigma^*}{R}{\Sigma}$--comodule $\Sigma$ (Theorem
\ref{comatrixconst}) is nothing but the original structure of
right $\coring{C}$--comodule over $\Sigma$. We have then a
commutative diagram of functors
\begin{equation}\label{CANSIGMA}
\xymatrix{\rcomod{\rcomatrix{R}{\Sigma}} \ar^{\mathsf{CAN_r}}[rr]
& & \rcomod{\coring{C}}\\
& \rmod{R} \ar^{\dostensor{-}{R}{\Sigma}}[lu]
\ar_{\dostensor{-}{R}{\Sigma}}[ru]&}
\end{equation}
On the other hand, the action of $\mathsf{CAN}_l$ on the left
$\dostensor{\Sigma^*}{R}{\Sigma}$--comodule $\Sigma^{\dag} =
\dostensor{\Sigma^*}{R}{R}$ (Theorem \ref{comatrixconst}) endows
it with a structure of left $\coring{C}$--comodule
$\lambda_{\Sigma^{\dag}}^{\coring{C}} : \Sigma^{\dag} \rightarrow
\dostensor{\coring{C}}{A}{\Sigma^{\dag}}$. Explicitly, it reads
\begin{equation*}
\lambda_{\Sigma^{\dag}}^{\coring{C}}(\dostensor{\phi}{R}{r}) =
\trestensor{\phi(e_{s[0]})e_{s[1]}}{A}{\varphi_s}{R}{r^s},
\end{equation*}
because, by definition, $\lambda_{\Sigma^{\dag}}^{\coring{C}} =
(\dostensor{\mathsf{can}}{A}{\Sigma^{\dag}})\lambda_{\Sigma^{\dag}}$.
Clearly, $\lambda_{\Sigma^{\dag}}^{\coring{C}}$ is right
$R$--linear and, thus, $\Sigma^{\dag} \in
\bicomod{\coring{C}}{R}$. If we assume ${}_RR$ to be flat (i.e.,
the functor $\dostensor{-}{R}{R}:\rmd{R} \rightarrow \mathsf{Ab}$
exact), this bicomodule gives a functor $- \cotensor{\coring{C}}
\Sigma^{\dag} : \rcomod{\coring{C}} \rightarrow \rmod{R}$. Recall
from Proposition \ref{homisotensor} that we have a functorial
isomorphism $\alpha : \dostensor{\hom{A}{\Sigma}{-}}{R}{R} \simeq
\dostensor{-}{A}{\Sigma^{\dag}}$.

\medskip

The following theorem establishes a basic adjunction connecting
comodules and firm modules.

\begin{theorem}\label{homcotensor}
If ${}_RR$ is flat, then the natural isomorphism $\alpha :
\dostensor{\hom{A}{\Sigma}{-}}{R}{R} \simeq
\dostensor{-}{A}{\Sigma^{\dag}}$ induces a natural isomorphism
$\dostensor{\hom{\coring{C}}{\Sigma}{-}}{R}{R} \simeq -
\cotensor{\coring{C}} \Sigma^{\dag}$. In particular, we have an
adjoint pair
\begin{equation}\label{adj5}
\xymatrix{\rmod{R} \ar@<0.5ex>^-{\dostensor{-}{R}{\Sigma}}[rrr] &
& & \rcomod{\coring{C}}, \ar@<0.5ex>^-{- \cotensor{\coring{C}}
\Sigma^{\dag}}[lll] & \dostensor{-}{R}{\Sigma} \dashv -
\cotensor{\coring{C}} \Sigma^{\dag}}.
\end{equation}
\end{theorem}
\begin{proof}
Consider the diagram of equalizers
\begin{equation*}
\xymatrix{\dostensor{\hom{\coring{C}}{\Sigma}{M}}{R}{R} \ar[d] & &
M
\cotensor{\coring{C}} \Sigma^{\dag} \ar[d] \\
\dostensor{\hom{A}{\Sigma}{M}}{R}{R}
\ar@<-1ex>[d]_{\dostensor{\beta''}{R}{R}}
\ar@<1ex>[d]^{\dostensor{\delta''}{R}{R}} \ar[rr]^-{\alpha_M}& &
\dostensor{M}{A}{\Sigma^{\dag}}
\ar@<-1ex>[d]_{\dostensor{\rho_M}{A}{\Sigma^{\dag}}}
\ar@<1ex>[d]^{\dostensor{M}{A}{\lambda_{\Sigma^{\dag}}}}
\\
 \dostensor{\hom{A}{\Sigma}{\dostensor{M}{A}{\coring{C}}}}{R}{R}
 \ar[rr]^-{\alpha_{\dostensor{M}{A}{\coring{C}}}}& & \trestensor{M}{A}{\coring{C}}{A}{\Sigma^{\dag}}, }
\end{equation*}
where $\beta''$ and $\delta''$ appeared in diagram
\eqref{equalizer3}. In order to deduce that $\alpha_{M}$ induces
an isomorphism $\alpha_M :
\dostensor{\hom{\coring{C}}{\Sigma}{M}}{R}{R} \simeq M
\cotensor{\coring{C}} \Sigma^{\dag}$ we need to check the
following two identities:\\
(a) $\alpha_{\dostensor{M}{A}{\coring{C}}} \circ
(\dostensor{\beta''}{R}{R}) =
(\dostensor{\rho_{\Sigma}}{A}{\Sigma^{\dag}})\circ \alpha_M$\\
(b) $\alpha_{\dostensor{M}{A}{\coring{C}}} \circ
(\dostensor{\delta''}{R}{R}) =
(\dostensor{M}{A}{\lambda_{\Sigma^{\dag}}}) \circ \alpha_M $.\\
To check (a), pick $\dostensor{h}{R}{r} \in
\dostensor{\hom{A}{\Sigma}{M}}{R}{R}$ and compute as follows.
\begin{equation*}
\begin{array}{rcl}
\alpha_{\dostensor{M}{A}{\coring{C}}}
(\dostensor{\beta''}{R}{R})(\dostensor{h}{R}{r}) & = &
\alpha_{\dostensor{M}{A}{\coring{C}}}(\dostensor{\rho_{\Sigma}h}{R}{r})
\\
& = & \trestensor{\rho_{\Sigma}h(e_s)}{A}{\varphi_s}{R}{r^s} \\
& = & \fourtensor{h(e_s)_{[0]}}{A}{h(e_s)_{[1]}}{A}{\varphi_s}{R}{r^s} \\
& = &
(\dostensor{\rho_{\Sigma}}{A}{\Sigma^{\dag}})(\trestensor{h(e_s)}{A}{\varphi_s}{R}{r^s})
\\
& = & (\dostensor{\rho_{\Sigma}}{A}{\Sigma^{\dag}})
\alpha_M(\dostensor{h}{R}{r}).
\end{array}
\end{equation*}
We shall deduce (b) by proving
$\alpha_{\dostensor{M}{A}{\coring{C}}} \circ
(\dostensor{\delta''}{R}{R}) \circ \alpha_M^{-1} =
\dostensor{M}{A}{\lambda_{\Sigma^{\dag}}}$ as follows: given
$\trestensor{m}{A}{\phi}{R}{r} \in
\dostensor{M}{A}{\Sigma^{\dag}}$, write $\dostensor{h}{R}{r} =
\alpha_M^{-1}(\trestensor{m}{A}{\phi}{R}{r})$, which means (see
the proof of Proposition \ref{homisotensor}) that $h(u) =
m\phi(u)$ for every $u \in \Sigma$. Therefore,
\begin{equation*}
\begin{array}{rcl}
(\alpha_{\dostensor{M}{A}{\coring{C}}} \circ
(\dostensor{\delta''}{R}{R}) \circ \alpha_M^{-1})
(\trestensor{m}{A}{\phi}{R}{r})& = &
\alpha_{\dostensor{M}{A}{\coring{C}}}
(\dostensor{\delta''}{R}{R})(\dostensor{h}{R}{r}) \\
 & = & \fourtensor{h(e_{s[0]})}{A}{e_{s[1]}}{A}{\varphi_s}{R}{r^s} \\
 & = & \fourtensor{m\phi(e_{s[0]})}{A}{e_{s[1]}}{A}{\varphi_s}{R}{r^s}
 \\
 & = & \fourtensor{m}{A}{\phi(e_{s[0]})e_{s[1]}}{A}{\varphi_s}{R}{r^s}
 \\
 & = & (\dostensor{M}{A}{\lambda_{\Sigma^{\dag}}})(\trestensor{m}{A}{\phi}{R}{r}).
\end{array}
\end{equation*}
This finishes the proof.
\end{proof}

The adjunction isomorphism of the adjoint pair \eqref{adj5} is
given by the restriction to
$\hom{\coring{C}}{\dostensor{N}{R}{\Sigma}}{M}$ of
\eqref{isoadj6}. This allows to deduce that the counit of the
adjoint pair \eqref{adj5} is given by the composite map
\begin{equation}\label{counit5}
\xymatrix{\chi_M : \dostensor{(M \cotensor{\coring{C}}
\Sigma^{\dag})}{R}{\Sigma}
\ar^-{\dostensor{eq_{M,\Sigma^{\dag}}}{R}{\Sigma}}[rr] & &
\trestensor{M}{A}{\Sigma^{\dag}}{R}{\Sigma} \ar^-{\delta_M}[r] &
M, }
\end{equation}
where $\delta_M$ is defined in \eqref{counit6}, and
$eq_{M,\Sigma^{\dag}} : M \cotensor{\coring{C}} \Sigma^{\dag}
\rightarrow \dostensor{M}{A}{\Sigma^{\dag}}$ is the equalizer of
the maps $\dostensor{\rho_M}{A}{\Sigma^{\dag}}$ and
$\dostensor{M}{A}{\lambda_{\Sigma^{\dag}}^{\coring{C}}}$.

\begin{proposition}\label{candos}
The following diagram is commutative.
\[
\xymatrix{ & & \Sigma^{\dag} \tensor{R} \Sigma \tensor{A}
\coring{C}
\ar^-{\dostensor{\pi_{A}}{A}{\coring{C}}}[rrd] & & \\
           \Sigma^{\dag} \tensor{R} \Sigma
           \ar[rru]^-{\Sigma^{\dag} \tensor{R}\rho_\Sigma}
           \ar_-{\lambda_{\Sigma^{\dag}} \tensor{R}\Sigma}[drr] \ar@{-->}^-{\mathsf{can}^{\dag}}[rrrr]& & & &\coring{C}\\
           & & \coring{C} \tensor{A} \Sigma^{\dag} \tensor{R} \Sigma \ar_-{\delta_{\coring{C}}}[urr] & }
\]
The resulting map
$$\mathsf{can}^{\dag}:\Sigma^{\dag} \tensor{R} \Sigma \rightarrow
\coring{C},$$ is such that the homomorphism of $A$--corings
$\mathsf{can}: \dostensor{\Sigma^*}{R}{\Sigma} \rightarrow
\coring{C}$ can be expressed as $\mathsf{can} =
\mathsf{can}^{\dag} \circ (\dostensor{\Sigma^*}{R}{d_{\Sigma}})$.
Moreover, $\mathsf{can}^{\dag}$ is a homomorphism of $A$--corings.
\end{proposition}
\begin{proof}
Pick $\trestensor{\phi}{R}{r}{A}{x} \in
\dostensor{\Sigma^{\dag}}{R}{\Sigma}$. Taking that $x = sx^{s}$
into account we get
\begin{equation*}
\begin{array}{rcl}
(\dostensor{\pi_A}{A}{\coring{C}})(\dostensor{\Sigma^{\dag}}{R}{\rho_{\Sigma}})(\trestensor{\phi}{R}{r}{A}{x})
& = &
(\dostensor{\pi_A}{A}{\coring{C}})(\dostensor{\Sigma^{\dag}}{R}{\rho_{\Sigma}})(\trestensor{\phi}{R}{r}{A}{sx^s})
\\
& = &
(\dostensor{\pi_A}{A}{\coring{C}})(\dostensor{\Sigma^{\dag}}{R}{\rho_{\Sigma}})(\trestensor{\phi}{R}{r}{A}{e_s\varphi_s(x^s)})
\\
& = &
(\dostensor{\pi_A}{A}{\coring{C}})(\fourtensor{\phi}{R}{r}{R}{e_{s[0]}}{A}{e_{s[1]}\varphi_s(x^s)})
\\
& = & \phi(re_{s[0]})e_{s[1]}\varphi_s(x^s),
\end{array}
\end{equation*}
and, since $\trestensor{\phi r}{R}{ts^t}{R}{x^s} =
\trestensor{\phi r}{R}{st}{R}{(x^s)^t} = \trestensor{\phi
r}{R}{s}{R}{x^s}$, we have
\begin{equation*}
\begin{array}{rcl}
\delta_{\coring{C}}(\dostensor{\lambda_{\Sigma^{\dag}}}{R}{\Sigma})(\trestensor{\phi}{R}{r}{R}{x})
& = &
\delta_{\coring{C}}(\dostensor{\lambda_{\Sigma^{\dag}}}{R}{\Sigma})(\trestensor{\phi}{R}{rs}{R}{x^s})
\\
& = &
\delta_{\coring{C}}(\trestensor{\mathsf{can}}{A}{\Sigma^{\dag}}{R}{\Sigma})(\phi
r\otimes_Re_s\otimes_A\varphi_s\otimes_Rt\otimes_R(x^s)^t) \\
 & = & \phi(re_{s[0]})e_{s[1]}\varphi_s(x^s).
\end{array}
\end{equation*}
The equality $\mathsf{can} = \mathsf{can}^{\dag} \circ
(\dostensor{\Sigma^*}{R}{d_{\Sigma}})$ follows without difficulty
from \eqref{cancompli}. Since
$\dostensor{\Sigma^*}{R}{d_{\Sigma}}$ is the inverse of the
isomorphism of $A$--corings given in Proposition \ref{comparar},
we get that $\mathsf{can}^{\dag}$ is a homomorphism of
$A$--corings.
\end{proof}

\begin{remark}\label{piC}
The canonical map $\mathsf{can}$ and the value of the counit
$$\pi_{\coring{C}} : \hom{\coring{C}}{\Sigma}{\coring{C}}
\tensor{R} R \tensor{R} \Sigma \rightarrow \coring{C}$$ at
$\coring{C}$ are tightly related. In fact, $\mathsf{can} =
\pi_{\coring{C}} \circ \daleth$, where $\daleth$ denotes the
composition of the isomorphisms
\begin{equation*}
\xymatrix{\rcomatrix{R}{\Sigma} \ar^-{\Sigma^* \tensor{R}
d_{\Sigma}}[rr] & & \Sigma^* \tensor{R} R \tensor{R} \Sigma
\ar^{((- \tensor{A} \coring{C}) \circ \rho_{\Sigma})\tensor{R} R
\tensor{R} \Sigma}[rrrr] & & & &
\hom{\coring{C}}{\Sigma}{\coring{C}}} \tensor{R} R \tensor{R}
\Sigma
\end{equation*}
\end{remark}

The following definition generalizes the notion of a Galois
comodule from the case of finite comatrix corings
\cite{ElKaoutit/Gomez:2003} to our general framework.

\begin{definition}\label{Galoisdef}
The comodule $\Sigma$ is said to be
$R$--$\coring{C}$-\emph{Galois} if $\mathsf{can}$ (or,
equivalently, $\mathsf{can}^{\dag}$) is an isomorphism. We say
then that $(\coring{C},{}_R\Sigma)$ is Galois, or even that
$\coring{C}$ is a Galois coring, when $\Sigma$ and $R$ are clear
from the context.
\end{definition}

From now on in this section, assume that ${}_RR$ is flat. We will
show that a Galois comodule $\Sigma$ allows to reconstruct some
comodules from the category $\rmod{R}$ by using the functor $-
\tensor{R} \Sigma$. To do this, consider the following diagram
($\Sigma$ is not yet assumed to be Galois)
\[
\xymatrix{  (N \cotensor{\coring{C}} \Sigma^{\dag}) \tensor{R}
\Sigma \ar^-{eq_{N,\Sigma^{\dag}} \tensor{R} \Sigma}[rrrr]
\ar^-{eq_{N,\Sigma^{\dag}}\tensor{R} \Sigma}[dd]
\ar@/_5pc/_-{\chi_N}[dddd]& & & & N \tensor{A} \Sigma^{\dag}
\tensor{R} \Sigma
\ar_-{N \tensor{A} \lambda_{\Sigma^{\dag}} \tensor{R} \Sigma}[dd] \ar@/^5pc/[dddd]^{N\tensor{A}\mathsf{can}^{\dag}}\\
& & & & \\
 N \tensor{A} \Sigma^{\dag} \tensor{R} \Sigma
\ar^{\rho_N \tensor{A} \Sigma^{\dag} \tensor{R} \Sigma}[rrrr]
\ar^-{\delta_N}[dd] & & & & N \tensor{A} \coring{C} \tensor{A}
\Sigma^{\dag} \tensor{R}
\Sigma \ar_-{\delta_{\dostensor{N}{A}{\coring{C}}}}[dd]\\
& & & & \\
N \ar_-{\rho_N}[rrrr]& & & & N \tensor{A} \coring{C}}
\]
which is commutative by Proposition \ref{candos}. Define a natural
transformation
\[
\xymatrix{\Psi : (- \cotensor{\coring{C}} \Sigma^{\dag})
\tensor{R} \Sigma \ar[r] & - \cotensor{\coring{C}} (\Sigma^{\dag}
\tensor{R} \Sigma)},
\]
making commute the diagrams (for $N \in \rcomod{\coring{C}}$)
\begin{equation}\label{hacha}
\xymatrix{(N \cotensor{\coring{C}} \Sigma^{\dag}) \tensor{R}
\Sigma \ar[ddrr]^{eq_{N,\Sigma^{\dag}}\tensor{R}\Sigma} \ar@{-->}[dd]^-{\Psi_N}& &  & \\
& & & \\
 N \cotensor{\coring{C}} (\Sigma^{\dag} \tensor{R} \Sigma)
\ar[rr]_{eq_{N,\Sigma^{\dag} \tensor{R} \Sigma}} & & N \tensor{A}
\Sigma^{\dag} \tensor{R} \Sigma \ar[r] \ar@<1ex>[r] & N \tensor{A}
\coring{C} \tensor{A} \Sigma^{\dag} \tensor{A} \Sigma}
\end{equation}
We thus get the following diagram
\begin{equation}\label{cinema}
\xymatrix{& & & N \tensor{A} \Sigma^{\dag} \tensor{R} \Sigma
\ar[ddd]^-{N \tensor{A} \mathsf{can}^{\dag}}
\\
& & & \\ (N \cotensor{\coring{C}} \Sigma^{\dag}) \tensor{R} \Sigma
\ar[rr]^-{\Psi_N} \ar[uurrr]^-{eq_{N,\Sigma^{\dag}} \tensor{R}
\Sigma} \ar[ddd]^-{\chi_N}& & N \cotensor{\coring{C}}
(\Sigma^{\dag} \tensor{R} \Sigma) \ar[ruu]_-{eq_{N,\Sigma^{\dag}
\tensor{R}
\Sigma}} \ar[ddd]_{N \cotensor{\coring{C}} \mathsf{can}^{\dag}}& \\
& & & N \tensor{A} \coring{C} \\
& & & \\
 N \ar_-{\sim}[rr] \ar'[rru]^-{\rho_N}[rrruu]& & N
\cotensor{\coring{C}} \coring{C} \ar_-{eq_{N,\coring{C}}}[uur]& }
\end{equation}
which turns out to be commutative.

\medskip

From the foregoing discussion we get the following
characterization of those comodules which come from modules.

\begin{theorem}\label{clave}
The following are equivalent for $N \in \rcomod{\coring{C}}$
\begin{enumerate}[(i)]
\item The counit $\zeta_N$ gives an isomorphism
$\dostensor{\hom{\coring{C}}{\Sigma}{N}}{R}{\Sigma} \cong N$;
\item the counit $\pi_N$ gives an isomorphism
$\trestensor{\hom{\coring{C}}{\Sigma}{N}}{R}{R}{R}{\Sigma} \cong
N$; \item the counit $\chi_N$ gives an isomorphism $(N
\cotensor{\coring{C}} \Sigma^{\dag}) \tensor{R} \Sigma \cong N.$
\end{enumerate}
If $\Sigma$ is an $R-\coring{C}$--Galois comodule, then the
foregoing statements are equivalent to each one of the following.
\begin{enumerate}
\item[(iv)] $\Psi_N$ gives an isomorphism $(N
\cotensor{\coring{C}} \Sigma^{\dag}) \tensor{R} \Sigma \cong N
\cotensor{\coring{C}} (\Sigma^{\dag} \tensor{R} \Sigma);$
\item[(v)] $- \tensor{R} \Sigma$ preserves the equalizer
$eq_{N,\Sigma^{\dag}}$.
\end{enumerate}
\end{theorem}
\begin{proof}
$(i) \Leftrightarrow (ii)$ follows from the commutative diagram
\[
\xymatrix{\dostensor{\hom{\coring{C}}{\Sigma}{N}}{R}{\Sigma}
\ar^-{\zeta_N}[rr] & & N \\
\trestensor{\hom{\coring{C}}{\Sigma}{N}}{R}{R}{R}{\Sigma}
\ar^-{\cong}[u] \ar_-{\pi_N}[urr] & & },
\]
where the vertical isomorphism comes from the isomorphism
$\dostensor{R}{R}{\Sigma} \cong \Sigma$. \\
$(ii) \Leftrightarrow (iii)$ follows from Theorem
\ref{homcotensor}.\\
$(iii) \Leftrightarrow (iv)$ is a consequence of diagram
\eqref{cinema}.\\
$(iv) \Leftrightarrow (v)$ is easily deduced from diagram
\eqref{hacha}.
\end{proof}

By combining Theorem \ref{clave} and Gabriel-Popescu's Theorem
\cite{Gabriel/Popescu:1964}, we obtain the following
generalization and improvement of \cite[Theorem
2.1.(1)]{Brzezinski:unp2004} and \cite[Theorem
3.8]{Caenepeel/DeGroot/Vercruysse:unp2005}. Observe that, as we
explain later, the use of the Gabriel-Popescu Theorem is not the
only possibility to deduce such a generalization. In fact, Theorem
\ref{flatdescent} can be understood as two different statements
that generalize the aforementioned results. Recall that the
category of comodules $\rcomod{\coring{C}}$ is a Grothendieck
category with exact forgetful functor $\rcomod{\coring{C}}
\rightarrow \rmod{A}$ if and only if ${}_A\coring{C}$ is flat
\cite[Proposition 1.2]{ElKaoutit/Gomez/Lobillo:2004}.

Let $T = \rend{\coring{C}}{\Sigma}$ be the endomorphism ring of
the right $\coring{C}$--comodule $\Sigma$. Since $\Sigma$ is an
$R-\coring{C}$--bicomodule, every element $r$ in the firm ring $R$
gives rise by left multiplication to an element of $T$. This
defines a homomorphism of rings $R \rightarrow T$.

\begin{theorem}\label{flatdescent}
The following statements are equivalent.
\begin{enumerate}[(i)]
\item \label{ff2} $\hom{\coring{C}}{\Sigma}{-} :
\rcomod{\coring{C}} \rightarrow \rmd{R}$ is a full and faithful
functor; \item $\dostensor{\hom{\coring{C}}{\Sigma}{-}}{R}{R} :
\rcomod{\coring{C}} \rightarrow \rmod{R}$ is a full and faithful
functor;  \item $-\cotensor{\coring{C}} \Sigma^{\dag} :
\rcomod{\coring{C}} \rightarrow \rmod{R}$ is a full and faithful
functor; \item $\Sigma$ is a Galois right $\coring{C}$--comodule
and the canonical map $\Psi_N$ gives an isomorphism $(N
\cotensor{\coring{C}} \Sigma^{\dag}) \tensor{R} \Sigma \cong N
\cotensor{\coring{C}} (\Sigma^{\dag} \tensor{R} \Sigma)$ for every
right $\coring{C}$--comodule $N$; \item \label{equalizer} $\Sigma$
is a Galois right $\coring{C}$--comodule and
$\dostensor{-}{R}{\Sigma}$ preserves the equalizer
$eq_{N,\Sigma^{\dag}}$ for every right $\coring{C}$--comodule $N$.

\medskip
Assume that $\coring{C}$ is flat as a left $A$--module, and
consider the following statements.  \item \label{sigmaplano}
$\Sigma$ is a Galois right $\coring{C}$--comodule and a flat left
$R$--module;
\item \label{generador} $\Sigma$ is a generator in
the category $\rcomod{\coring{C}}$; \item \label{ffT}
$\hom{\coring{C}}{\Sigma}{-} : \rcomod{\coring{C}} \rightarrow
{\Mm_T}$ is a full and faithful functor.
\end{enumerate}
Then (vii) and (viii) are equivalent, (vi) implies (v), and (i)
implies (vii). If, in addition, ${}_RT$ is flat or $R$ is a left
ideal of $T$, then (vii) implies (vi) and, therefore, the eight
statements are equivalent.
\end{theorem}
\begin{proof}
The equivalence between the first five statements is a direct
consequence of Theorem \ref{clave}, in view of Remark \ref{piC}.\\
$(vii) \Leftrightarrow (viii)$ is a consequence of
Gabriel-Popescu's Theorem \cite[Cap. III, Teorem\u a 9.1.(2)]{Nastasescu:1976}.\\
$(vi) \Rightarrow (v)$ is obvious. \\
$(i) \Rightarrow (vii)$ is clear. \\
$(viii) \Rightarrow (vii)$ The main problem is to prove that
${}_R\Sigma$ is flat. We will prove this here in the case ${}_RT$
flat. The other case, namely for $R$ a left ideal of $T$ will be
proved below in Proposition \ref{plano1}. By Gabriel-Popescu's
Theorem \cite[Cap. III, Teorem\u a 9.1.(3)]{Nastasescu:1976}, the
functor $\dostensor{-}{T}{\Sigma} : \rmod{T} \rightarrow
\rcomod{\coring{C}}$ is exact. Since the forgetful functor
$\rcomod{\coring{C}} \rightarrow \rmod{A}$ is exact because
${}_A\coring{C}$ is flat, we get that ${}_T\Sigma$ is flat. This
implies, being ${}_RT$ flat, that ${}_R\Sigma$ is flat. Now,
Theorem \ref{clave} gives that $\chi_{\coring{C}}$ is an
isomorphism. From \eqref{cinema} we get that $\coring{C}
\cotensor{\coring{C}} \mathsf{can}^{\dag}$ is an isomorphism and,
hence, $\mathsf{can}^{\dag}$ is an isomorphism. Therefore,
$\Sigma$ is Galois.
\end{proof}

Apart from the equivalence between $(vi)$ and $(viii)$ in Theorem
\ref{flatdescent}, Gabriel-Popescu's Theorem is used in $(viii)
\Rightarrow (vii)$ just to deduce that ${}_R\Sigma$ is flat from
the fact that $\Sigma_{\coring{C}}$ is a generator, in the case
that ${}_RT$ is flat. For the other  case, when $R$ is a left
ideal of $T$, we need to extend some technical facts on flat
modules from unital rings to firm rings. This is made, for the
convenience of the reader, in the Appendix. With these results at
hand, we proceed as follows.

Let $R \rightarrow T$ be a homomorphism of rings from a firm ring
$R$ to a unital ring $T$. We have a functor $\dostensor{T}{R}{-} :
\lmod{R} \rightarrow \lmod{T}$, and, for every $N \in \lmod{R}$,
an homomorphism of left $R$--modules $v_N : N \rightarrow
\dostensor{T}{R}{N}$ defined as $v_N(n) = \dostensor{1}{R}{n}$.

\begin{lemma}\label{RT}
The following statements are equivalent.
\begin{enumerate}[(i)]
\item $v_R : R \rightarrow \dostensor{T}{R}{R}$ is an isomorphism;
\item $R$ is a left ideal of $T$; \item $v_N : N \rightarrow
\dostensor{T}{R}{N}$ is an isomorphism for every $N \in \lmod{R}$.
\end{enumerate}
\end{lemma}
\begin{proof}
$(i) \Rightarrow (ii)$ For $t \in T$ and $r \in R$ put $s =
v_R^{-1}(\dostensor{t}{R}{r}) \in R$. Then $\dostensor{1}{R}{s} =
\dostensor{t}{R}{r}$ which implies, by multiplying, that $s =
tr$.\\ $(ii) \Rightarrow (i)$ Define $f : \dostensor{T}{R}{R}
\rightarrow R$ by $f(\dostensor{t}{R}{r}) = tr$. Obviously, $f
\circ v_R = R$. For the other composition,
\[
(v_R \circ f)(\dostensor{t}{R}{r}) = v_R(tr) = v_R(tsr^s) =
tsv_R(r^s) = \dostensor{ts}{R}{r^s} = \dostensor{t}{R}{sr^s} =
\dostensor{t}{R}{r}
\]
which shows that $f = v_R^{-1}$. \\
$(i) \Rightarrow (iii)$ With ${}_RN$ firm we have
\[
\dostensor{T}{R}{N} \cong \trestensor{T}{R}{R}{R}{N} \cong
\dostensor{R}{R}{N} \cong N,
\]
isomorphism which is explicitly given by
\begin{equation}\label{haceTmodulo}
\dostensor{t}{R}{n} \mapsto \trestensor{t}{R}{r}{R}{n^r} \mapsto
\dostensor{tr}{R}{n^r} \mapsto (tr)n^r
\end{equation}
This isomorphism is the inverse of $v_N$ since, for $n \in N$,
\[
n \mapsto \dostensor{1}{R}{n} \mapsto \trestensor{1}{R}{r}{R}{n^r}
\mapsto \dostensor{r}{R}{n^r} \mapsto rn^r = n
\]
$(iii) \Rightarrow (i)$ is obvious.
\end{proof}

\begin{lemma}\label{leftflat}
If $R$ is a left ideal of $T$ then every firm left $R$--module $N$
is a left $T$--module with the action $tn = (tr)n^r$ for $t \in T$
and $n \in N$. Moreover, we have an isomorphism
$\dostensor{M}{T}{N} \cong \dostensor{M}{R}{N}$ for every $M \in
\rmod{T}$.
\end{lemma}
\begin{proof}
The isomorphism $v_N$ transfers the structure of left $T$--module
from $\dostensor{T}{R}{N}$ to $N$ accordingly to
\eqref{haceTmodulo}. The isomorphism $\dostensor{M}{T}{N} \cong
\dostensor{M}{R}{N}$ is then obtained, with the help of Lemma
\ref{RT}, as the composite
\[
\dostensor{M}{T}{N} \cong \dostensor{M}{T}{(\dostensor{T}{R}{N})}
\cong \dostensor{(\dostensor{M}{T}{T})}{R}{N} \cong
\dostensor{M}{R}{N}
\]
\end{proof}

\begin{remark}
If $R$ is a left ideal of $T = \rend{\coring{C}}{\Sigma}$ then, by
Lemma \ref{leftflat}, we have an isomorphism
$\rcomatrix{T}{\Sigma} \cong \rcomatrix{R}{\Sigma}$ which turns
out to be of $A$--bimodules. Therefore, the structure of comatrix
$A$--coring of $\rcomatrix{R}{\Sigma}$ can be transferred to
$\rcomatrix{T}{\Sigma}$, in analogy to the case of infinite
comatrix corings \cite{ElKaoutit/Gomez:2004}.
\end{remark}

We are ready to formulate the alternative proof for the flatness
of $\Sigma$. It is modelled after \cite[19.9]{Wisbauer:1991}. We
assume that ${}_A\coring{C}$ is flat.

\begin{proposition}\label{plano1}
If $\Sigma$ is a generator of $\rcomod{\coring{C}}$ and $R$ is a
left ideal of $T = \rend{\coring{C}}{\Sigma}$, then ${}_R\Sigma$
is flat.
\end{proposition}
\begin{proof}
By Proposition \ref{flat4} in the Appendix it suffices to show
that for any right ideal of finite type $J=f_{1}R+\cdots+f_kR$ of
$R$, the map $\mu_J:J\otimes_R\Sigma\to J\Sigma$, $\mu_J(g\otimes
u)=g(u)$ is injective. Let us consider the surjective map
$\phi:\Sigma^n\to J\Sigma,f(u_{1},\ldots,u_n)=\sum_if_i(u_i)$. We
put $K=\Ker\phi$, then $K\in\Mm^\coring{C}$, since
${_A\coring{C}}$ is flat. Moreover, we have an exact sequence
\[
\xymatrix{ 0\ar[r] &\hom{\coring{C}}{\Sigma}{K} \ar[r]^\alpha &
\hom{\coring{C}}{\Sigma}{\Sigma^n} \ar[r]^-\beta & J \ar[r] & 0. }
\]
Tensoring by $\Sigma$ over $R$, we obtain the following commutative diagram with exact rows:
\[
\xymatrix{ &\hom{\coring{C}}{\Sigma}{K}\otimes_R\Sigma
\ar[r]^{\alpha\otimes\Sigma} \ar[d]^{{\rm ev}_K} &
\hom{\coring{C}}{\Sigma}{\Sigma^n}\otimes_R\Sigma \ar[r]^-{\beta\otimes\Sigma} \ar[d]^{{\rm ev}_{\Sigma^n}}& J\otimes_R\Sigma \ar[r] \ar[d]^{\mu_J} & 0\\
0 \ar[r] & K \ar[r] & \Sigma^n \ar[r]^\phi & J\Sigma \ar[r] &0
}
\]
The map ${\rm ev}_K$ is surjective, this is a translation of the
generator property of $\Sigma$ (see, e.g., \cite[Lemma
3.7]{Caenepeel/DeGroot/Vercruysse:unp2005}). Furthermore, ${\rm
ev}_{\Sigma^n}$ is an isomorphism, as one can see from the
following
$$\hom{\coring{C}}{\Sigma}{\Sigma^n}\otimes_R\Sigma\cong \hom{\coring{C}}{\Sigma}{\Sigma}^n\otimes_R\Sigma\cong
 T^n\otimes_R\Sigma\cong T^n\otimes_RR\otimes_R\Sigma\cong R^n\otimes_R\Sigma\cong \Sigma^n,$$
(here we used Lemma \ref{RT}). After this, it follows from the
properties of the diagram that $\mu_J$ is injective.
\end{proof}

We have then finished the proof of Theorem \ref{flatdescent}. The
following consequence allow to deduce (see Theorem
\ref{descensofielmenteplano}) a generalization of the faithfully
flat descent theorem for comodules as formulated in
\cite{ElKaoutit/Gomez:2003}. For the notion of a quotient category
we refer to \cite{Gabriel:1962}.

\begin{corollary}\label{Gabriel}
If $\coring{C}$ is flat as left $A$ module and if the equivalent conditions of Theorem \ref{flatdescent} hold,
then $\hom{\coring{C}}{\Sigma}{-}\otimes_RR$ induces an equivalence of categories between $\Mm^\coring{C}$ and
$\Mm_R/\ker(-\otimes_R\Sigma)$.
\end{corollary}

\begin{proof}
This is an immediate consequence of \cite[Proposition
III.5]{Gabriel:1962}.
\end{proof}

The following is our main theorem, which generalizes \cite[Theorem
3.2]{ElKaoutit/Gomez:2003} and, as we will show in Section
\ref{comparison}, \cite[Theorem 5.9]{ElKaoutit/Gomez:2004}.

\begin{theorem}\label{descensofielmenteplano} Let $\coring{C}$ be an $A$-coring that is flat as left
$A$-module. With notations as before, the following statements are
equivalent
\begin{enumerate}[(i)]
\item $\Sigma$ is a Galois right $\coring{C}$-comodule and a
faithfully flat left $R$-module; \label{Sigmaffl} \item
$-\otimes_R\Sigma:\Mm_R\to\Mm^\coring{C}$ is an equivalence of
categories; \label{equiv} \item $\Sigma$ is a Galois right
$\coring{C}$-comodule and
$-\otimes_R\Sigma:\Mm_R\to\Mm^{\Sigma^*\otimes_R\Sigma}$ is an
equivalence of categories; \label{equivcom} \item $\Sigma$ is a
generator of $\rcomod{\coring{C}}$ such that the functor
$\dostensor{-}{R}{\Sigma} : \rmod{R} \rightarrow
\rcomod{\coring{C}}$ is full and faithful; \label{gene1} \item
$\Sigma$ is a generator of $\rcomod{\coring{C}}$ such that the
functor $\dostensor{-}{R}{\Sigma} : \rmod{R} \rightarrow
\rcomod{\coring{C}}$ is faithful and $R$ is a left ideal of $T$.
\label{gene2}
\end{enumerate}
\end{theorem}
\begin{proof}
\textit{(\ref{Sigmaffl})} $\Rightarrow$ \textit{(\ref{equiv})} If
$\Sigma$ is faithfully flat as a left $R$--module, then the kernel
of $\dostensor{-}{R}{\Sigma}$ is trivial. By Corollary
\ref{Gabriel}, $\dostensor{-}{R}{\Sigma}$ is an equivalence of
categories.\\
$(\ref{equiv}) \Rightarrow (\ref{equivcom})$ If
$\dostensor{-}{R}{\Sigma} : \rmod{R} \rightarrow
\rcomod{\coring{C}}$ is an equivalence, then its pseudo-inverse
$\dostensor{\hom{\coring{C}}{\Sigma}{-}}{R}{R}$ is an equivalence
and, by Theorem \ref{flatdescent}, $\Sigma$ is Galois, i.e.,
$\mathsf{can} : \rcomatrix{R}{\Sigma} \rightarrow \coring{C}$ is
an isomorphism of corings. Then the induced functor
$\mathsf{CAN_r} : \rcomod{\rcomatrix{R}{\Sigma}} \rightarrow
\rcomod{\coring{C}}$ is an equivalence of categories. This
implies, in view of diagram \eqref{CANSIGMA} that
$\dostensor{-}{R}{\Sigma} : \rmod{R} \rightarrow
\rcomod{\coring{C}}$ is an equivalence. \\
$(\ref{equivcom}) \Rightarrow (\ref{equiv})$ follows from diagram
\eqref{CANSIGMA}.\\
$(\ref{equiv}) \Rightarrow (\ref{gene1})$ Since $R$ is a generator
for $\rmod{R}$, we get that $\Sigma \cong
\dostensor{R}{R}{\Sigma}$ is a generator of $\rcomod{\coring{C}}$.\\
$(\ref{gene1}) \Rightarrow (\ref{gene2})$ The unit of the
adjunction \eqref{adj4} is an isomorphism. This gives an
isomorphism
$$R\cong \dostensor{\hom{\coring{C}}{\Sigma}{R\otimes_R\Sigma}}{R}{R}
\cong\dostensor{\hom{\coring{C}}{\Sigma}{\Sigma}}{R}{R}=\dostensor{T}{R}{R},$$
which entails that every element in $\dostensor{T}{R}{R}$ is of
the form $\dostensor{\mu_R\nu_R(r^s)}{R}{s}$ for some $r \in R$
(notations as in \eqref{unit3}). In this way, the image of the
multiplication map $\dostensor{T}{R}{R} \rightarrow T$ is $R$,
that is, $R$ is a left ideal of $T$.\\
$(\ref{gene2}) \Rightarrow (\ref{Sigmaffl})$ The forgetful functor
$\rcomod{\coring{C}} \rightarrow \rmod{A}$ is faithful, which
implies, under the hypotheses, that $\dostensor{-}{R}{\Sigma} :
\rmod{R} \rightarrow \rmod{A}$ is faithful. By Theorem
\ref{flatdescent}, ${}_R\Sigma$ is flat and it is a Galois
comodule.
\end{proof}

Let us return to the case where $\coring{C} =
\rcomatrix{R}{\Sigma}$. That $(\rcomatrix{R}{\Sigma},{}_R\Sigma)$
is Galois is not a surprise, even thought that it is not
completely evident.

\begin{lemma}\label{canid}
The homomorphism of $A$--corings $\mathsf{can} :
\rcomatrix{R}{\Sigma} \rightarrow \rcomatrix{R}{\Sigma}$ is the
identity. Therefore, $\Sigma$ is a Galois
$\rcomatrix{R}{\Sigma}$--comodule.
\end{lemma}
\begin{proof}
By using  \eqref{Raccion} and the definition of the coaction
$\rho_{\Sigma}$, we have
\begin{equation*}
\begin{array}{rcl}
\mathsf{can}(\dostensor{\phi}{R}{x}) & = & \phi(e_r)\varphi_r
 \tensor{R} x^r \\
 & = & \phi r \tensor{R} x^r \\
 & = & \phi \tensor{R} rx^r \\
 & = & \phi \tensor{R} x
\end{array}
\end{equation*}
\end{proof}

In view of the Lemma \ref{canid}, Theorem
\ref{descensofielmenteplano} has a relevant corollary.

\begin{corollary}\label{ffdescent}
Assume that the comatrix coring $\rcomatrix{R}{\Sigma}$ is flat as
a left $A$--module. The following statements are equivalent.
\begin{enumerate}[(i)]
\item $\Sigma$ is a faithfully flat left $R$-module; \item
$-\otimes_R\Sigma:\Mm_R\to\Mm^{\Sigma^*\otimes_R\Sigma}$ is an
equivalence of categories.
\end{enumerate}
\end{corollary}

\begin{remark}
The statements of this section admit a reformulation in terms of
the coring $\Sigma^{\dag} \tensor{R} \Sigma$ defined in
\eqref{deltadag}. This claim is supported by Proposition
\ref{comparar} and Proposition \ref{candos}, which provide the
commutative diagram of homomorphisms of $A$--corings, with $f$ an
isomorphism:
\begin{equation}\label{canplus}
\xymatrix{\rcomatrix{R}{\Sigma} \ar^{\mathsf{can}}[dr] & \\
& \coring{C} \\
\Sigma^{\dag} \tensor{R} \Sigma \ar^{f}[uu]
\ar_{\mathsf{can}^{\dag}}[ru] & }
\end{equation}
\end{remark}

\begin{example}
Suppose $\Sigma$ is a $B-\coring{C}$-bicomodule, such that
the following condition holds. Every finite subset of $\Sigma$ is
contained in a $B-\coring{C}$--subcomodule, say $P$, of $\Sigma$
such that $P_A$ is finitely generated and projective and a direct
summand of ${}_B\Sigma_{\coring{C}}$. If $\pi : \Sigma \rightarrow
P$ is the corresponding projection and $\{(e_i ,e_i^*)\} \subseteq
P \times P^*$ is a dual basis for $P_A$, then $e = e_i \tensor{A}
e_i^*\pi \in\dostensor{\Sigma}{A}{\Sigma^*}$ is an idempotent such
that $e$ is $B$-central, $e$ is an endomorphism of the
$\coring{C}$--comodule $\Sigma$, and $e$ acts as a left unit on
$P$. This situation is more restrictive than the one in Example
\ref{locprojcom}. In case $B=\mathbb{Z}$ and $\coring{C}=A$, it
means exactly that $\Sigma_A$ is locally projective in the sense
of \'Anh and M\'arki \cite{AM}. As in Example \ref{locprojcom} we
can construct a comatrix coring $\dostensor{\Sigma^*}{R}{\Sigma}$,
where $R$ is a subring of $\dostensor{\Sigma}{A}{\Sigma^*}$
generated by a bunch of idempotent dual bases for subcomodules of
$\Sigma$, indexed by a set of generators of $\Sigma$. Each of this
idempotents is an endomorphism of $\Sigma_{\coring{C}}$, and,
therefore, $\Sigma$ becomes an $R-\coring{C}$--bicomodule. This
implies we can apply the results of this section to this
situation. We can as well take for $R$ the left ideal $T$
generated by the idempotents, which would make this extra
condition in Theorem 4.9 superfluous. One can prove that the
conditions in this example are equivalent with taking $\Sigma$
equal to a split direct limit of $\coring{C}$-comodules that are
finitely generated and projective. In this way, our results cover
the ones of \cite{CDV2005}.
\end{example}

\section{Infinite Comatrix Corings}\label{comparison}

Next, we will show that one of the constructions of the infinite
comatrix corings of \cite{ElKaoutit/Gomez:2004}, and the most
relevant results of \cite[Section 5]{ElKaoutit/Gomez:2004} can be
deduced from the theory developed here. With this purpose, let us
start by reconsidering Example \ref{infcomatrix}. There, $\Sigma =
\oplus_{P \in \cat{P}} P$ for a set $\cat{P}$ of finitely
generated and projective right $A$--modules. Write $\Sigma^{\ddag}
= \oplus_{P \in \cat{P}}P^*$. We see that $\Sigma^{\dag} =
\Sigma^* \tensor{R} R$ is isomorphic, as an $A-R$--bimodule, to
$\Sigma^{\ddag}$. In fact, for every $P \in \cat{P}$, we have an
isomorphism of left $A$--modules $P^* \cong \Sigma^{\dag}
\tensor{R} R$ that sends $\phi$ onto $\phi \pi_P \tensor{R} u_P$.
This set of isomorphisms defines an isomorphism of
$A-R$--bimodules $h : \Sigma^{\dag} \rightarrow \Sigma^{\ddag}$
which, for its part, induces the isomorphism of $A$--bimodules
\begin{equation}\label{chicoring}
h \tensor{R} \Sigma : \Sigma^{\dag} \tensor{R} \Sigma \rightarrow
\Sigma^{\ddag} \tensor{R} \Sigma
\end{equation}
The structure of $A$--coring of $\Sigma^{\dag} \tensor{R} \Sigma$
given in \eqref{Deltadag2} is obviously transferred via $h
\tensor{R} \Sigma$ to $\Sigma^{\ddag} \tensor{R} \Sigma$. In
$\Sigma^{\ddag} \tensor{R} \Sigma$ we have that
$\dostensor{\phi}{R}{x} = 0$ for every $\phi \in P^*$ and $x \in
Q$, whenever $P \neq Q$ (see \cite[page
2032]{ElKaoutit/Gomez:2004}). Therefore, those elements obtained
as sums of $\phi \tensor{R} x$, where $\phi \in P^*$ and $x \in
P$, exhaust $\Sigma^{\ddag} \tensor{R} \Sigma$ when $P$ runs
$\cat{P}$. The comultiplication $\Delta^{\ddag}$ of
$\Sigma^{\ddag} \tensor{R} \Sigma$ is then determined by
\begin{equation*}
\Delta^{\ddag}(\dostensor{\phi}{R}{x}) = \sum
\fourtensor{\phi}{R}{e_{\alpha_P}}{A}{e_{\alpha_P}^*}{R}{x} \qquad
( \phi \in P^*, x \in P, P \in \cat{P}),
\end{equation*}
in concordance with the structure given in \cite[Proposition
5.2]{ElKaoutit/Gomez:2004}. The isomorphism of corings
\eqref{chicoring} induces an $R-\Sigma^{\ddag} \tensor{R}
\Sigma$--bicomodule (resp. a $\Sigma^{\ddag} \tensor{R} \Sigma
-R$--bicomodule) structure on $\Sigma$ (resp. $\Sigma^{\ddag}$).
In resume, in virtue of Proposition \ref{comparar} and the
foregoing remarks, the statements of Section \ref{Galoiscomodules}
have some consequences on the coring $\Sigma^{\dag}\tensor{R}
\Sigma$ (or equivalently, on $\Sigma^{\ddag} \tensor{R} \Sigma$)
which complete the theory developed in \cite[Section
5]{ElKaoutit/Gomez:2004}. We will give the most interesting of
them.

\medskip

Let $\cat{P}$ be a set of right comodules over an $A$--coring
$\coring{C}$, and define $\Sigma = \oplus_{P \in \cat{P}} P$.
Consider the ring $T = \rend{\coring{C}}{\Sigma}$ and its (non
unital) subring $R = \oplus_{P,Q \in \cat{P}} u_PTu_Q$. When each
$P \in \cat{P}$ is finitely generated and projective as a right
$A$--module, then we have the $A$--coring $\Sigma^{\dag}
\tensor{R} \Sigma$ and the homomorphism of $A$--corings
$\mathsf{can}^{\dag} : \Sigma^{\dag} \tensor{R} \Sigma \rightarrow
\coring{C}$. On the other hand, it is not difficult to prove that
${}_RR$ is flat and that $R$ is a left ideal of $T$. We are then
in the hypotheses of Theorem \ref{flatdescent}.

\begin{theorem}\label{predescensoplano}
Let $\Sigma = \oplus_{P \in \cat{P}}P$ a right comodule over an
$A$--coring $\coring{C}$, where $\cat{P}$ is a set of right
$\coring{C}$--comodules that are finitely generated and projective
as right $A$--modules. Consider $R = \bigoplus_{P,Q \in \cat{P}}
u_PTu_Q$, where $u_P$ is the idempotent of $T =
\rend{\coring{C}}{\coring{C}}$ ``attached'' to $P$. The following
statements are equivalent.
\begin{enumerate}[(i)]
\item \label{ff2b} $\hom{\coring{C}}{\Sigma}{-} :
\rcomod{\coring{C}} \rightarrow \rmd{R}$ is a full and faithful
functor; \item $\dostensor{\hom{\coring{C}}{\Sigma}{-}}{R}{R} :
\rcomod{\coring{C}} \rightarrow \rmod{R}$ is a full and faithful
functor;  \item $-\cotensor{\coring{C}} \Sigma^{\dag} :
\rcomod{\coring{C}} \rightarrow \rmod{R}$ is a full and faithful
functor; \item $\Sigma$ is a Galois right $\coring{C}$--comodule
and the canonical map $\Psi_N$ gives an isomorphism $(N
\cotensor{\coring{C}} \Sigma^{\dag}) \tensor{R} \Sigma \cong N
\cotensor{\coring{C}} (\Sigma^{\dag} \tensor{R} \Sigma)$ for every
right $\coring{C}$--comodule $N$; \item \label{equalizerb}
$\Sigma$ is a Galois right $\coring{C}$--comodule and
$\dostensor{-}{R}{\Sigma}$ preserves the equalizer
$eq_{N,\Sigma^{\dag}}$ for every right $\coring{C}$--comodule $N$.

\medskip
Assume that $\coring{C}$ is flat as a left $A$--module. Then the
foregoing statements are equivalent to each of the following.

\item \label{sigmaplanob} $\Sigma$ is a Galois right
$\coring{C}$--comodule and a flat left $R$--module;

\item \label{generadorb} $\cat{P}$ is a set of generators in the
category $\rcomod{\coring{C}}$; \item \label{ffTb}
$\hom{\coring{C}}{\Sigma}{-} : \rcomod{\coring{C}} \rightarrow
{\Mm_T}$ is a full and faithful functor.
\end{enumerate}
\end{theorem}

The main result of \cite[Section 5]{ElKaoutit/Gomez:2004} is a
consequence of our general theory.

\begin{theorem}\label{descensoplano}\cite[Theorem
5.7]{ElKaoutit/Gomez:2004} Let $\cat{P}$ be a set of right
comodules over an $A$--coring $\coring{C}$, and consider $\Sigma =
\oplus_{P \in \cat{P}}P$. Let $R = \oplus_{P \in \cat{P}}
u_PTu_Q$, where $T = \rend{\coring{C}}{\Sigma}$, as before. The
following statements are equivalent
\begin{enumerate}[(i)]
\item ${}_A{\coring{C}}$ is flat and $\cat{P}$ is a generating set
of small projectives for $\rcomod{\coring{C}}$; \item
${}_A\coring{C}$ is flat, every $P \in \cat{P}$ is finitely and
generated and projective as a right $A$--module,
$\mathsf{can}^{\dag} : \Sigma^{\dag} \tensor{R} \Sigma \rightarrow
\coring{C}$ is an isomorphism, and $- \tensor{R} \Sigma :
\rcomod{\Sigma^{\dag} \tensor{R} \Sigma} : \rmod{R} \rightarrow
\rcomod{\coring{C}}$ is an equivalence of categories; \item every
$P \in \cat{P}$ is finitely generated and projective as a right
$A$--module, $\mathsf{can}^{\dag} : \Sigma^{\dag} \tensor{R}
\Sigma \rightarrow \coring{C}$ is an isomorphism, and ${}_R\Sigma$
is faithfully flat; \item ${}_A\coring{C}$ is flat and $-
\tensor{R} \Sigma : \rmod{R} \rightarrow \rcomod{\coring{C}}$ is
an equivalence of categories.
\end{enumerate}
\end{theorem}
\begin{proof}
This is a consequence of Theorem \ref{descensofielmenteplano},
after the following observation: from $(i)$ it follows that every
$P \in \cat{P}$ is finitely generated and projective as a right
$A$--module, as it is proved in \cite[Theorem
5.7]{ElKaoutit/Gomez:2004}.
\end{proof}

\section{Comatrix corings applied to rationality
properties}\label{rationality}

Let $\coring{C}$ be an $A$-coring with counity
$\varepsilon_{\coring{C}}$. The left dual
${}^*\coring{C}=\hom{A}{{}_A\coring{C}}{A}$ has the structure of
an $A$-ring with multiplication
\begin{equation}
\label{mult*c} f*g(c)=g(c_{(1)}f(c_{(2)}))
\end{equation}
and unit map
$$\iota:A\to {}^*\coring{C},\qquad \iota(a)(c)=\varepsilon_{\coring{C}}(c)a,$$
which induces an $A$-bimodule structure on ${}^*\coring{C}$ given
by
\begin{equation}
\label{multA*c} (fa)(c)=f(c)a,\qquad (af)(c)=f(ca).
\end{equation}
Every right $\coring{C}$-comodule $M$ is has also a right
${}^*\coring{C}$-module structure given by
$$m\cdot f=m_{[0]}f(m_{[1]}),$$
for all $m\in M$ and $f\in {}^*\coring{C}$. If $\coring{C}$ is
locally projective as a left $A$-module, then we can partially
converse this last property and study rational
${}^*\coring{C}$-modules \cite{A, CVW2}. For any
$M\in\rmod{{}^*\coring{C}}$, we denote by $M^{\rm rat}$ the
rational submodule of $M$. This is given by
$$M^{\rm rat}=\{m\in M~|~\exists m_i\in M, c_i\in \coring{C}, \textrm{ s.t. } m\cdot f=m_if(c_i), \forall f\in{}^*\coring{C}\},$$
in this case $M^{\rm rat}$ is a right $\coring{C}$-comodule with
coaction given by
\begin{equation}
\label{coactionrat} \rho(m)=m_{[0]}\otimes_Am_{[1]}\quad
\textrm{iff}\quad m\cdot f=m_{[0]}f(m_{[1]}), \textrm{ for all }
f\in {}^*\coring{C}.
\end{equation}
By $\cat{R}\rmod{^*\coring{C}}$ we denote the full subcategory of
$\rmod{{}^*\coring{C}}$ consisting of all rational submodules.
There exists  an isomorphism of categories
\begin{equation}
\label{isocat}
\rcomod{\coring{C}}\cong \cat{R}\rmod{{}^*\coring{C}}
\end{equation}
Recall from \cite{A, Caenepeel/DeGroot/Vercruysse:unp2005, CVW2}
that ${}^*\coring{C}^{\rm rat}$ is a ring with right local units
if and only if it is a dense subset in the finite topology of
${}^*\coring{C}$. In this situation, ${}^*\coring{C}^{\rm rat}$
has also right local units on every right $\coring{C}$-comodule
$M$. It is clear that $\coring{C}\in\Mm_{{}^*\coring{C}^{\rm
rat}}$ and ${}^*\coring{C}^{\rm rat}\in{{}^*{\coring{C}_{\rm
rat}}\Mm}$ and by the presence of local units, both are firm
${}^*\coring{C}^{\rm rat}$-modules. In this situation the
isomorphism of categories \eqref{isocat} can be extended with an
equivalence
\begin{equation}
\Mm^\coring{C}\cong \cat{R}\rmod{{}^*\coring{C}} \simeq
\Mm_{({}^*\coring{C})^{\rm rat}},
\end{equation}
where $\rmod{{}^*\coring{C}^{\rm rat}}$ denotes the category of
firm right ${}^*\coring{C}^{\rm rat}$-modules, which means in this
situation that ${{}^*\coring{C}^{\rm rat}}$ acts with right local
units on them. The aim of this section is to show that the
composed equivalence $\rcomod{\coring{C}} \simeq
\rmod{{}^*\coring{C}^{\rm rat}}$ can be obtained by use of
 Galois comodules.

Consider as before the $A$-coring $\coring{C}$ which is locally
projective as left $A$-module.

\begin{lemma}
\label{lemma1} The following equation holds for all $c\in
\coring{C}$ and $f\in {}^*\coring{C}^{\rm rat}$
\begin{equation}
\label{question} c_{(1)}f(c_{(2)})=f_{[0]}(c)f_{[1]}
\end{equation}
\end{lemma}

\begin{proof}
Let $\{e_i,f_i\}$ be a local final dual basis for both
$c_{(1)}f(c_{(2)})$  and $f_{[0]}(c)f_{[1]}$. Then
\[
\begin{array}{rcll}
f_{[0]}(c)f_{[1]}&=&e_if_i(f_{[0]}(c)f_{[1]})\\
&=&e_if_{[0]}(c)f_i(f_{[1]})&\textrm{$f_i$ left $A$-linear}\\
&=&e_i(f_{[0]}f_i(f_{[1]}))(c)&\eqref{multA*c}\\
&=&e_i(f*f_i)(c)&\eqref{coactionrat}\\
&=&e_if_i(c_{(1)}f(c_{(2)}))&\eqref{mult*c}\\
&=&c_{(1)}f(c_{(2)})
\end{array}
\]
\end{proof}

\begin{theorem}
The right coaction on ${}^*\coring{C}^{\rm rat}$, given by
$$\rho:{}^*\coring{C}^{\rm rat}\to {}^*\coring{C}^{\rm rat}\otimes_A\coring{C},\qquad f\mapsto f_{[0]}\otimes f_{[1]},$$
is a homomorphism of rings.
\end{theorem}

\begin{proof}
The multiplication on ${}^*\coring{C}^{\rm rat}$ is induced by
\eqref{mult*c}. The multiplication on ${}^*\coring{C}^{\rm
rat}\otimes_A\coring{C}$ is defined as usual by
$$(f\otimes_A c)\cdot (g\otimes_A d)=f\otimes_A g(c)d.$$
We can compute
\begin{equation}
\label{ringmorph1}
\begin{array}{rcll}
\rho(f*g)&=&(f*g)_{[0]}\otimes_A(f*g)_{[1]}\\
&=&f*g_{[0]}\otimes_Ag_{[1]}&\textrm{\cite[Prop. 5.1]{CVW2}}
\end{array}
\end{equation}
and
\begin{equation}
\label{ringmorph2}
\begin{array}{rcl}
\rho(f)\cdot\rho(g)&=&(f_{[0]}\otimes_Af_{[1]})\cdot(g_{[0]}\otimes_Ag_{[1]})\\
&=&f_{[0]}g_{[0]}(f_{[1]})\otimes_Ag_{[1]}
\end{array}
\end{equation}
Furthermore, $f*g_{[0]}(c)=g_{[0]}(c_{(1)}f(c_{2}))$ and
$$(f_{[0]}g_{[0]}(f_{[1]}))(c)=f_{[0]}(c)g_{[0]}(f_{[1]})=g_{[0]}(f_{[0]}(c)f_{[1]}).$$ Finally, using \eqref{question}, we find that the
right hand sides of \eqref{ringmorph1} and \eqref{ringmorph2} are
equal and thus $\rho$ is indeed a homomorphism of rings.
\end{proof}

From now on, we will restrict ourselves  to the situation where
${}^*\coring{C}^{\rm rat}$ is a ring with right local units.

Now we can apply the results of Section \ref{comatrixsection} and
construct a comatrix coring
$\coring{C}\otimes_{{}^*\coring{C}^{\rm rat}}{{}^*\coring{C}^{\rm
rat}}$. The comultiplication and counit are given respectively by
\[
\begin{array}{rrcl}
\Delta : & \coring{C}\otimes_R{{}^*\coring{C}^{\rm rat}} &\to& \coring{C}\otimes_R{{}^*\coring{C}^{\rm rat}}\otimes_A\coring{C}\otimes_R{{}^*\coring{C}^{\rm rat}}\\
&c\otimes_R f & \mapsto & c\otimes_R f_{[0]}\otimes_A f_{[1]}\otimes_R e \\
\varepsilon : &\coring{C}\otimes_R{{}^*\coring{C}^{\rm rat}} &\to & A\\
 & c\otimes_R f & \mapsto & f(c)
\end{array}
\]
where $e$ is a right local unit for $f$ and we denote
$R={{}^*\coring{C}^{\rm rat}}$ from now on, if we concern only its
ring structure.

The canonical map is given by the map of \ref{lemma1} :
$$\can :  \coring{C}\otimes_R{{}^*\coring{C}^{\rm rat}} \to \coring{C}, \quad \can(c\otimes_R f)=f_{[0]}(c)f_{[1]}\quad(=c_{(1)}f(c_{(2)})\quad),$$
and has an inverse that reads
$$\can^{-1}(c)=c\otimes_Re,$$
where $e\in R$ is a local unit for $c$ (i.e. $c\cdot e= c$).
Remark that this isomorphism is nothing else the firmness
isomorphism $\coring{C}\cong \coring{C}\otimes_RR$.

We will now apply the results of Section \ref{Galoiscomodules}.
The Galois comodule ${}^*\coring{C}^{\rm rat}$ induces a pair of
adjoint functors given by
\begin{equation}
\label{adjoint} \xymatrix{ \Mm_R
\ar@<0.5ex>[rrr]^-{-\otimes_R{}^*\coring{C}^{\rm rat}} &&&
\Mm^\coring{C}
\ar@<0.5ex>[lll]^{\Hom^\coring{C}({}^*\coring{C}^{\rm
rat},-)\otimes_RR} }
\end{equation}

The general theory makes use of the module
$\Sigma^\dagger=\Sigma^*\otimes_RR$. In our situation, we have the
following

\begin{lemma}
\label{dagger} With notations and conventions as before, the
following holds,
$${}^*\coring{C}^{{\rm rat}\dagger}=\Hom_A({}^*\coring{C}^{\rm rat},A)\otimes_RR\cong \coring{C}.$$
\end{lemma}

\begin{proof}
Take $\varphi\otimes_R r\in\Hom_A({}^*\coring{C}^{\rm
rat},A)\otimes_RR$ then we define
$\alpha(\varphi\otimes_Rr)=\varphi(r_{[0]})r_{[1]}$, and
conversely we define for all $c\in\coring{C}$,
$\beta(c)=\psi_c\otimes_r e$ where $e\in R$ is a right local unit
for $c$ and $\psi_c$ is defined by $\psi_c(f)=f(c)$ for all
$f\in{}^*\coring{C}^{\rm rat}$. Then we can check
\[
\begin{array}{rcll}
\alpha\circ\beta(c)&=&e_{[0]}(c)e_{[1]}\\
&=&c_{(1)}e(c_{(2)})&\eqref{question}\\
&=&c\cdot e= c
\end{array}
\]
and
\begin{equation}
\label{betaalpha}
\begin{array}{rcl}
\beta\circ\alpha(\varphi\otimes_R
r)&=&\psi_{\varphi(r_{[0]})r_{[1]}}\otimes_R e
\end{array}
\end{equation}
where $e$ is a right local unit for $\varphi(r_{[0]})r_{[1]}$ and
thus also for $r$. Let us first compute
\[
\begin{array}{rcl}
\psi_{\varphi(r_{[0]})r_{[1]}}(f)&=&f(\varphi(r_{[0]})r_{[1]})\\
&=&\varphi(r_{[0]})f(r_{[1]})=\varphi(r_{[0]}f(r_{[1]}))\\
&=&\varphi(r*f)=(\varphi\cdot r)(f)
\end{array}
\]
We can now go on with \eqref{betaalpha},
\[
\begin{array}{rcl}
\beta\circ\alpha(\varphi\otimes_R r)&=&\varphi\cdot r\otimes_R e\\
&=&\varphi\otimes_Rr\cdot e=\varphi\otimes_R r.
\end{array}
\]
\end{proof}

Since the existence of local units in $R$ implies that $_RR$ is
flat, by Theorem 4.4 and the previous lemma,
$$\hom{\coring{C}}{{}^*\coring{C}^{\rm rat}}{-}\otimes_RR\simeq -\Box_\coring{C}\coring{C},$$
however, let us give the isomorphisms of this equivalence for
completeness sake, they are very similar to the ones in the proof
of Lemma \ref{dagger}. Take $M\in M^\coring{C}$ then
$$\alpha:\Hom_\coring{C}({}^*\coring{C}^{\rm rat},M)\otimes_RR\to M\Box_\coring{C}\coring{C}\cong M,$$
is an isomorphism where for all
$\varphi\otimes_Rr\in\Hom_\coring{C}({}^*\coring{C}^{\rm
rat},M)\otimes_RR$, $m\in M$, $\alpha(\varphi\otimes
r)=\varphi(r)$ and $\alpha^{-1}(m)=\psi_m\otimes e$, where
$\psi_m(f)=m\cdot f$ for all $f\in R$ and $e$ is a right local
unit for $m$.

It follows now that the pair of adjoint functors induced by the
Galois comodule ${}^*\coring{C}^{\rm rat}$ is can be written as
\[
\xymatrix{ \Mm_R \ar@<0.5ex>[rrr]^-{-\otimes_RR} &&&
\Mm^\coring{C} \ar@<0.5ex>[lll]^{-\Box_\coring{C}\coring{C}} }
\]
Since $M\otimes_RR\cong M$ for every firm right $R$-module $M$ and
$N\Box_\coring{C}\coring{C}\cong N$ for every right
$\coring{C}$-comodule $N$, it is clear that these functors
constitute the same equivalence of categories as mentioned in the
beginning of this section.

\section{Appendix: Flat modules over firm rings}

In this section we will always consider $R$ to be a firm ring. Our
aim is to generalize some well-known facts about flat and
injective modules over unital rings to the situation of firm
rings. Our treatment is an adaptation of the one given in
\cite[Chapter I, Section 10]{Stenstrom:1975}.

Recall that $M$ is called an injective right $R$-module if and
only if $M$ is an injective object in the category $\rmd{R}$,
which translates to the the fact that for every injective morphism
$i:N\to N'$ of right $R$-modules, the induced map
$\Hom_R(N',M)\to\Hom_R(N,M)$ must be surjective.

\begin{proposition}
\label{flat1} A right $R$-module $E$ is injective if and only if
$\Hom_R(R,E)\to\Hom(I,E)$ is surjective for every right ideal $I$
of $R$.
\end{proposition}

\begin{proof}
In one direction the statement is trivial. For the other way, we
start from a monomorphism $\alpha:L\to M$ and any morphism
$\varphi:L\to E$. Consider the set
$$\Mm=\{\varphi':L'\to E~|~L\subset L'\subset M, \textrm  { s.t.  $\varphi'$ extends $\varphi$}\}.$$
We can equip $\Mm$ with a partial order, defining
$\varphi'\le\varphi''$ if and only if $\varphi''$ is a further
extension.

Take a totally ordered $\Tt\subset\Mm$ and put
$\tilde{L}=\sum_{L'\in \Tt} L'$, and define
$\tilde{\varphi}:\tilde{L}\to E$ by the condition
$\tilde{\varphi}(l)=\varphi'(l)$ for $l\in L'\in\Tt$. In this way
we find an upper bound for every totally ordered $\Tt\subset\Mm$
and by the Zorn's Lemma, this implies the existence of a maximal
element $\varphi_{0}:L_{0}\to E$ in $\Mm$. We have finished if we
can prove that $L_{0}=M$.

Suppose there exists a $x\in M$, $x\not\in L_{0}$, we show it is
possible to extend $\varphi_{0}$ to $\psi:L_{0}+xR\to E$, from
which the contradiction follows. Put $I=\{a\in R~|~xa\in L_{0}\}$,
then $I$ is a right ideal in $R$. Define $\beta:I\to E,
\beta(a)=\varphi_{0}(xa)$. By hypothesis, there exists a
$\varphi_R:R\to E$, extending $\beta$. For all $z\in L_{0}$ and
$a\in R$, we define $\psi(z+xa)=\varphi_{0}(z)+\varphi_R(a)$. To
see that $\psi$ is well defined, suppose $z+xa=z'+xa'$ then
$z-z'=x(a'-a)$ and we find $a'-a\in I$. By definition of
$\varphi_R$ we get then
$\varphi_{0}(z-z')=\varphi_{0}(x(a'-a))=\varphi_R(a'-a)$ and thus
$\psi(z+xa)=\varphi_{0}(z)+\varphi_R(a)=\varphi_{0}(z')+\varphi_R(a')=\psi(z'+xa')$.

Finally, $\psi$ is right $R$-linear and clearly it extends
$\varphi_{0}$, which ends the proof.
\end{proof}

\begin{proposition}
\label{flat2} Consider $F$ an $R-B$--bimodule, where $B$ is a
unital ring, and let $E$ be an injective cogenerator for
$\rmod{B}$. Then $F$ is flat as a left $R$-module if and only if
$\Hom_B(F,E)$ is injective as a right $R$-module
\end{proposition}

\begin{proof}
The proof is identical as in the classical case, since it makes no
use of units at all (see e.g. \cite{Stenstrom:1975}).
\end{proof}

\begin{corollary}
\label{flat3} A right $R$-module $F$ is flat if and only if
$\hat{F}=\Hom_\ZZ(F,\QQ/\ZZ)$ is injective as right $R$-module.
\end{corollary}

A right ideal $I$ of $R$ is said to be \emph{of finite type} if $I
= Ra_{1} + \cdots + Ra_n$ for some $a_{1}, \dots, a_n \in R$.

\begin{proposition}
\label{flat4} A firm left $R$-module $F$ is flat if and only if
the multiplication map $\mu_I:I\otimes_R F\to F$ is injective for
every ideal of finite type $I$ of $R$.
\end{proposition}

\begin{proof}
If $F$ is flat, the statement is clear, so we only prove the other
implication. Suppose that $\mu_I$ is injective for every right
ideal $I$ of finite type of $R$, then it is also the case for all
right ideals $J$ of $R$. Indeed, take $r_i\otimes y_i\in
J\otimes_RF$ and let $r_i = r_i^{s_i}s_i$ be given by the equality
$R^2 = R$. If we denote by $I = \sum r_i^{s_i} R$, which is of
finite type, then the following diagram is commutative and $\mu_I$
is injective,
\[
\xymatrix{
I\otimes_RF\ar[rr]^{\mu_I} \ar[dr]& & F\\
&J\otimes_RF\ar[ur]_{\mu_J} }
\]
It follows that $\mu_J(r_i\otimes
y_i)=0\Leftrightarrow\mu_I(r_i\otimes y_i)=0\Leftrightarrow
r_i\otimes y_i=0$, so $\mu_J$ is injective as well.

Now we will construct a commutative diagram
\[
\xymatrix{
\Hom_{R}(R,\hat{F})\ar[r]^\gamma \ar[d]_\cong^\alpha & \Hom_{R}(J,\hat{F}) \ar[d]^\cong_\beta \\
\hat{F} \ar[r]_\delta &\widehat{J\otimes_R F}, }
\]
where, for any $M$, we define $\hat{M}=\Hom_\ZZ(M,\QQ/\ZZ)$. For
all $\Phi\in\Hom_R(R,\hat{F})$ we define define
$\alpha(\Phi)(y)=\Phi(r)(y^r)$ for all $y\in F$. Conversely, we
define $\alpha^{-1}(\Psi)(r)(y)=\Psi(ry)$ for all $\Psi\in\hat{F},
r\in R, y\in F$. One easily computes
\[
\begin{array}{rcl}
\alpha\circ\alpha^{-1}(\Psi)(y)&=&\alpha^{-1}(\Psi)(r)(y^r)\\
&=&\Psi(ry^r)=\Psi(y)\\
\alpha^{-1}\circ\alpha(\Phi)(r)(y)&=&\alpha(\Phi)(ry)\\
&=&\Phi(s)((ry)^s)=\Phi(r)(y).
\end{array}
\]
In the last step, we used the $R$-linearity of $\Phi$ together
with $s\otimes_R(ry)^s=r\otimes_Ry$, which is based on the unique
expression of the element $ry$ mapped trough the isomorphism $d_M
: M\to R\otimes_RM$. Analogously, define
$\beta(\dostensor{a}{R}{y}) = \Phi(ar)(y^r)$ for every $\Phi \in
\hom{R}{J}{\hat{F}}$ and $\dostensor{a}{R}{y} \in J \tensor{R} F$,
which is shown to be an isomorphism analogously to the case of
$\alpha$.

Since $\mu_J:J\otimes F\to F$ is mono and $\QQ/\ZZ$ is injective,
$\delta$, which is the dual morphism of $\mu_J$ must be an
epimorphism and by the diagram, $\gamma$ is also epi. It follows
from Proposition \ref{flat1} that $\hat{F}$ is injective and thus
$F$ is flat by Corollary \ref{flat3}.
\end{proof}

\subsection*{Acknowledgement}
Joost Vercruysse would like to thank the Department of Algebra,
University of Granada for hospitality during his stay, where this
work initiated. This research has been partially supported by the
grant BFM2004-01406 from the Ministerio de Educaci{\'o}n y Ciencia of
Spain.

\providecommand{\bysame}{\leavevmode\hbox
to3em{\hrulefill}\thinspace}
\providecommand{\MR}{\relax\ifhmode\unskip\space\fi MR }
\providecommand{\MRhref}[2]{%
  \href{http://www.ams.org/mathscinet-getitem?mr=#1}{#2}
} \providecommand{\href}[2]{#2}

\end{document}